\documentclass{amsart}

\usepackage{amssymb}

\newtheorem{theorem}{Theorem}[section]
\newtheorem{corollary}[theorem]{Corollary}
\newtheorem{lemma}[theorem]{Lemma}
\newtheorem{proposition}[theorem]{Proposition}

\theoremstyle{definition}
\newtheorem{definition}[theorem]{Definition}
\newtheorem{example}[theorem]{Example}

\theoremstyle{remark}

\newcommand{\abs}[1]{\lvert#1\rvert}

\begin{document}

\title{Non-Commutative Metrics on Matrix State Spaces}

\author{wei wu}
\address{Department of Mathematics, East China Normal University, Shanghai 200062, P.R. China}
\email{wwu@math.ecnu.edu.cn}
\curraddr{Department of Mathematics, University of California, Berkeley, CA 94720-3840}
\email{wwu@math.berkeley.edu}

\begin{abstract}We use the theory of quantization to introduce non-commutative versions of
metric on state space and Lipschitz seminorm. We show that a lower
semicontinuous matrix Lipschitz seminorm is determined by their
matrix metrics on the matrix state spaces. A matrix metric comes
from a lower semicontinuous matrix Lip-norm if and only if it is
convex, midpoint balanced, and midpoint concave. The operator space of Lipschitz functions with a
matrix norm coming from a closed matrix Lip-norm is the operator
space dual of an operator space. They generalize Rieffel's results
to the quantized situation.
\end{abstract}

\maketitle

\section{Introduction}\label{sec:1}

In the non-commutative geometry,  the natural way to specify a
metric is by means of a suitable ``Lipschitz seminorm"\cite{co, con}. Given a triple $(\mathcal{A, H}, D)$, where $\mathcal
A$ is a unital $C^*$-algebra, $D$ is the generalized Dirac
operator on the Hilbert space $\mathcal H$ such that $(\mathcal H,
D)$ is an unbounded Fredholm module over $\mathcal A$. The set
$\mathcal{L(A)}$ of  Lipschitz elements of $\mathcal A$ consists
of those $a\in\mathcal A$ such that the commutator $[D, a]$ is a
bounded operator on $\mathcal H$. The Lipschitz seminorm, $L$, is
defined on $\mathcal{L(A)}$ just by the operator norm $L(a)=\|[D,
a]\|$. Using $L$, Connes defined a metric on the state space
$\mathcal{S(A)}$ of $\mathcal A$ by the formula
\[\rho_L(\varphi, \psi)=\sup\{|\varphi(a)-\psi(a)|:\, a\in\mathcal{L(A)}, L(a)\le 1\}.\]
This metric generalizes the Monge-Kantorovich metric on probability measures.

In \cite{ri}, Rieffel extended this idea to the order unit space and discussed the
relationship between the metrics on the state space and Lipschitz seminorms. A {\it matrix
order unit space} $(\mathcal V, 1)$ is a matrix ordered space $\mathcal V$ together with a
distinguished order unit $1$ satisfying the following conditions:
\begin{enumerate}
\item $\mathcal V^+$ is a proper cone with the order unit $1$;
\item each of the cones $M_n(\mathcal V)^+$ is Archimedean.
\end{enumerate}
The basic representation theorem of Kadison\cite{ka} says that any order unit space is
order isomorphic to a subspace of $C(X)$, the continuous functions on a compact
set $X$, closed under conjugation and containing the identity function. As the non-commutative
analogues of Kadison's {\it function systems}, {\it operator systems}, a self-adjoint linear
space of operators on a Hilbert space which contains the identity operator,  have been
characterized by Choi and Effros as matrix order unit spaces: every matrix order unit space is
completely order isomorphic to an operator system\cite{chef}. It is therefore only natural to
ask whether the metrics on the state space and Lipschitz seminorms generalize to this
non-commutative setting, and in this paper we give an affirmative answer to this question.

We begin with the explanation of terminology and notation. We
proceed by defining the non-commutative analogue of Lipschitz
seminorm and discuss some elementary properties of it. The next
step is to prove the non-commutative version of the fact that the
lower semicontinuous Lipschitz seminorms are determined by their
metrics on the state spaces. In section \ref{sec:5}, we define the
non-commutative analogue of the metric on compact convex sets, and
discuss the relation of matrix metrics on matrix state space and
matrix norms on dual operator space. Using the results of section
\ref{sec:5}, we prove, in section \ref{sec:6}, that as in the
commutative case, a matrix metric comes from a lower
semicontinuous matrix Lip-norm if and only if it is convex,
midpoint balanced, and midpoint concave.
We conclude in section \ref{sec:7} with an operator analogue of
the well-known result that the space of Lipschitz functions with a
norm coming from some Lipschitz seminorm is the dual of a certain
other Banach space.

\section{Terminology and notation}\label{sec:2}

All vector spaces are assumed to be complex throughout this paper.
Given a vector space $V$, we let $M_{m, n}(V)$ denote the matrix
space of all $m$ by $n$ matrices $v=[v_{ij}]$ with $v_{ij}\in V$,
and we set $M_n(V)=M_{n,n}(V)$. If $V=\mathbb C$,  we write
$M_{m,n}=M_{m,n} (\mathbb C)$ and $M_n=M_{n,n}(\mathbb C)$, which
means that we may identify $M_{m,n}(V)$ with the tensor product
$M_{m,n}\otimes V$. We identify $M_{m,n}$ with the normed space
$\mathcal B({\mathbb C}^n, {\mathbb C}^m)$. We use the standard
matrix multiplication and *-operation for compatible scalar
matrices, and $1_n$ for the identity matrix in $M_n$.

There are two natural operations on the matrix spaces. For $v\in
M_{m,n}(V)$ and $w\in M_{p,q} (V)$, the direct sum $v\oplus w\in
M_{m+p,n+q}(V)$ is defined by letting
\[v\oplus w=\left[\begin{array}{cc}v&0\\ 0&w\end{array}\right],\]
and if we are given $\alpha\in M_{m,p}$, $v\in M_{p,q}(V)$ and $\beta\in M_{q,n}$, the matrix
product  $\alpha v\beta\in M_{m,n}(V)$ is defined by
\[\alpha v\beta=\left[\sum_{k,l}\alpha_{ik}v_{kl}\beta_{lj}\right].\]

A *-{\it vector space} $V$ is a complex vector space together with
a conjugate linear mapping $v\longmapsto v^*$ such that
$v^{**}=v$. A complex vector space $V$ is said to be {\it matrix
ordered} if:
\begin{enumerate}
\item $V$ is a *-vector space;
\item each $M_n(V)$, $n\in\mathbb N$, is partially ordered;
\item  $\gamma^*M_n(V)^+\gamma\subseteq M_m(V)^+$ if $\gamma=[\gamma_{ij}]$ is any $n\times m$
matrix of complex numbers.
\end{enumerate}
A {\it matrix order unit space} $(\mathcal V, 1)$ is a matrix ordered space $\mathcal V$
together with a distinguished order unit $1$ satisfying the following conditions:
\begin{enumerate}
\item $\mathcal V^+$ is a proper cone with the order unit $1$;
\item each of the cones $M_n(\mathcal V)^+$ is Archimedean.
\end{enumerate}
Each matrix order unit space $(\mathcal V, 1)$ may be provided with the norm
\[\|v\|=\inf\left\{t\in\mathbb R:\, \left[\begin{array}{cc} t1&v\\ v^*&t1\end{array}
\right]\ge 0\right\}.\]
As in \cite{ri}, we will not assume that $\mathcal V$ is complete for the norm.

If $V$ and $W$ are *-vector spaces and $\varphi: V\longmapsto W$
is a linear mapping, we have a linear mapping $\varphi^*:
V\longmapsto W$ defined by $\varphi^*(v)={\varphi(v^*)}^*$.

Given vector spaces $V$ and $W$ and a linear mapping $\varphi:\,
V\longmapsto W$ and $n\in\mathbb N$, we have a corresponding
$\varphi_n:\, M_n(V)\longmapsto M_n(W)$ defined by
\[\varphi_n([v_{ij}])=[\varphi(v_{ij})].\]
If $V$ and $W$ are vector spaces in duality, then they determine
the matrix pairing
\[\ll\cdot, \cdot\gg:\, M_n(V)\times M_m(W)\longmapsto M_{nm},\]
where \[\ll[v_{ij}], [w_{kl}]\gg=\left[<v_{ij}, w_{kl}>\right]\]
for $[v_{ij}]\in M_n(V)$ and $[w_{kl}]\in M_m(W)$.

A {\it graded set} $\mathbf S=(S_n)$ is a sequence of sets
$S_n(n\in\mathbb N)$. If $V$ is a locally convex topological
vector space, then the canonical topology on $M_n(V)(n\in\mathbb
N)$ is that determined  by the natural linear isomorphism
$M_n(V)\cong V^{n^2}$, that is, the product topology. A graded set
$\mathbf S=(S_n)$ with $S_n\subseteq M_n(V)$ is {\it closed} if
that is the case for each set $S_n$ in the product topology in
$M_n(V)$. Given a vector space $V$, we say that a graded set
$\mathbf B=(B_n)$ with $B_n\subseteq M_n(V)$ is {\it absolutely
matrix convex} if for all $m, n\in\mathbb N$
\begin{enumerate}
\item $B_m\oplus B_n\subseteq B_{m+n}$;
\item $\alpha B_m\beta\subseteq B_n$ for any contractions $\alpha\in M_{n,m}$ and $\beta\in
M_{m,n}$.
\end{enumerate}
Let $V$ and $W$ be vector spaces in duality, and let $\mathbf
S=(S_n)$ be a graded set with $S_n\subseteq M_n(V)$. The {\it
absolute operator polar} ${\mathbf
S}^{\circledcirc}=(S_n^{\circledcirc})$ with
$S_n^\circledcirc\subseteq M_n(W)$, is defined by
$S_n^\circledcirc=\{w\in M_n(W):\,\|\ll v, w\gg\|\le 1\hbox{ for
all }v\in S_r, r\in\mathbb N\}$.

A {\it gauge} on a vector space $V$ is a function
$g:\,V\longmapsto [0, +\infty]$ such that
\begin{enumerate}
\item $g(v+w)\le g(v)+g(w)$;
\item $g(\alpha v)\le\abs{\alpha} g(v)$,
\end{enumerate}
for all $v,w\in V$ and $\alpha\in\mathbb C$. We say that a gauge $g$ is a {\it seminorm} on
$V$ if $g(v)<+\infty$ for all $v\in V$. Given an arbitrary vector space $V$, a {\it matrix
gauge} $\mathcal G=(g_n)$ on $V$ is a sequence of gauges
\[g_n:\, M_n(V)\longmapsto [0, +\infty]\]
such that
\begin{enumerate}
\item $g_{m+n}(v\oplus w)=\max\{g_m(v), g_n(w)\}$;
\item $g_n(\alpha v\beta)\le\|\alpha\| g_m(v)\|\beta\|$,
\end{enumerate}
for any $v\in M_m(V)$, $w\in M_n(V)$, $\alpha\in M_{n,m}$ and
$\beta\in M_{m,n}$. A matrix gauge $\mathcal G=(g_n)$ is a {\it
matrix seminorm} on $V$ if for any $n\in\mathbb N, g_n(v)<+\infty$
for all $v\in M_n(V)$. If each $g_n$ is a norm on $M_n(V)$, we say
that $\mathcal{G}$ is a {\it matrix norm}. An {\it operator space}
is a vector space together with a matrix norm on it. For a matrix
order unit space $(\mathcal V, 1)$, it is an operator space with
the matrix norm determined by the matrix order on it.

\section{Matrix Lipschitz seminorm}\label{sec:3}

We recall that a {\it Lipschitz seminorm} on an order unit space
$S$ is a seminorm on $S$ such that its null space is the scalar
multiples of the order unit\cite{ri}. When we have a triple
$(\mathcal{A, H}, D)$ as in the introduction, we obtain a sequence
of triples $(M_n({\mathcal A}), \mathcal H^n, D_n)$ and a sequence
of seminorms $L_n(a)=\|[D_n, a]\|$ on each $M_n(\mathcal{L(A)})$.
These seminorms are linked by the following fundamental relations:
\begin{enumerate}
\item the null space of each $L_n$ contains $M_n({\mathbb C1})$;
\item $L_{m+n}(v\oplus w)=\max\{L_m(v), L_n(w)\}$;
\item $L_n(\alpha v\beta)\le\|\alpha\| L_m(v)\|\beta\|$;
\item $L_m(v^*)=L_m(v)$,
\end{enumerate}
where $1$ is the identity of $\mathcal A$, $v\in
M_m(\mathcal{L(A)})$, $w\in M_n(\mathcal{L(A)})$, $\alpha \in
M_{n,m}$ and $\beta \in M_{m,n}$. Inspired by the observation, we
give our non-commutative analogue of a Lipschitz
seminorm\cite{wu}.

\begin{definition}\label{def:31} Given a matrix order unit space $(\mathcal V, 1)$, a {\it matrix
Lipschitz seminorm} $\mathcal L$ on $(\mathcal V, 1)$ is a sequence of seminorms
\[L_n:\,M_n(\mathcal V)\longmapsto [0, +\infty)\]
such that
\begin{enumerate}
\item the null space of each $L_n$ is $M_n({\mathbb C1})$;
\item $L_{m+n}(v\oplus w)=\max\{L_m(v), L_n(w)\}$;
\item $L_n(\alpha v\beta)\le\|\alpha\| L_m(v)\|\beta\|$;
\item $L_m(v^*)=L_m(v)$,
\end{enumerate}
for any $v\in M_m(\mathcal V)$, $w\in M_n(\mathcal V)$, $\alpha \in M_{n,m}$ and
$\beta \in M_{m,n}$.
\end{definition}

\begin{example}\label{ex:32}
Let $G$ be a compact group with identity element $e$. We suppose
that $G$ is equipped with a length function $l$, that is, a
continuous non-negative real-valued function on $G$ such that
\begin{enumerate}
\item $l(xy)\le l(x)+l(y)$,
\item $l(x^{-1})=l(x)$,
\item $l(x)=0$ if and only if $x=e$,
\end{enumerate}
for all $x,y\in G$.

Let $\alpha$ be an ergodic action of $G$ on a unital $C^*$-algebra
$\mathcal{A}$. For $a=[a_{ij}]\in M_n(\mathcal{A})$ and
$n\in\mathbb{N}$, we set
\[L_n(a)=\sup\left\{\frac{\|[\alpha_{x}(a_{ij})-a_{ij}]\|}{l(x)}: x\neq e\right\}.\]
We let $\mathcal{V}$ denote the set of Lipschitz elements of
$\mathcal{A}$ for $\alpha$ and $l$ with corresponding Lipschitz
seminorm $L_1$. Then $\mathcal{L}=(L_n)$ is a matrix Lipschitz
seminorm on $(\mathcal{V}, 1)$, where $1$ is the identity element
of $\mathcal{A}$.
\end{example}

If $\mathcal V$ is an operator space, we let
$\mathcal{CB}(\mathcal V, M_n)$ denote the vector space of
completely bounded linear mapping $\varphi:\,\mathcal V\longmapsto
M_n$, on which we place the completely bounded norm
$\|\cdot\|_{cb}$. By the isometric identification
\[M_n(\mathcal V^*)\cong\mathcal{CB}(\mathcal V, M_n),\, n\in\mathbb N,\]
the matrix norms on $\mathcal V^*$ determine the operator space
dual $\mathcal V^*$ of $\mathcal V$(see Lemma 2.1 in \cite{efrua} and Proposition 2.1 in \cite{ble}). This operator
space dual and the matrix pairing substitute the Banach space dual
and the scalar pairing in the transition from the commutative to
the non-commutative case.

Let $(\mathcal V, 1)$ be a matrix order unit space. The {\it
matrix state space} of $(\mathcal V, 1)$ is the collection
$\mathcal{CS}(\mathcal V)=(CS_n(\mathcal V))$ of {\it matrix
states}
\[CS_n(\mathcal V)=\{\varphi\in\mathcal{CB}(\mathcal V, M_n):\,\varphi \hbox{ is completely positive },\, \varphi
(1)=1_n\},\]\cite{wewi,far}. If ${\mathcal L}=(L_n)$ is a matrix
Lipschitz seminorm on $(\mathcal V, 1)$, then, for each
$n\in\mathbb{N}$ we can define a metric $D_{L_n}:\,CS_n(\mathcal
V)\times CS_n(\mathcal V)\longmapsto [0, +\infty]$ by
\[D_{L_n}(\varphi, \psi)=\sup\{\|\ll \varphi, a\gg -\ll\psi, a\gg\|:\, a\in M_r(\mathcal V),
L_r(a)\le 1, r\in{\mathbb N}\},\] where $\varphi , \psi  \in
CS_n(\mathcal V)$ and we should notice that it may take value
$+\infty$. Denote $\mathcal{D_L}=(D_{L_n})$. And in turn we can
define a sequence ${\mathcal L}_{{\mathcal D}_{\mathcal
L}}=(L_{D_{L_n}})$ of gauges on $(\mathcal V, 1)$ by
\[L_{D_{L_n}}(a)=\sup\left\{\frac{\|\ll \varphi, a\gg -\ll\psi, a\gg\|}{D_{L_r}(\varphi,
\psi)}:\,\varphi, \psi\in CS_r(\mathcal V), \varphi\neq \psi,
r\in{\mathbb N}\right\},\] for all $a\in M_n(\mathcal V)$.

\begin{proposition}\label{pro:33} Let ${\mathcal L}=(L_n)$ be a matrix Lipschitz seminorm on the matrix
order unit space $(\mathcal V, 1)$. Then we have
\[D_{L_n}(\varphi, \psi)=\sup\{\|\ll \varphi, a\gg -\ll\psi, a\gg\|:\, a=a^*\in M_r(\mathcal
V), L_r(a)\le 1, r\in{\mathbb N}\},\] for $\varphi, \psi\in
CS_n(\mathcal V)$.
\end{proposition}

\begin{proof} Assume that $D_{L_n}(\varphi, \psi)=c<+\infty$. For arbitrary $\epsilon >0$, there
exist $r\in\mathbb N$ and $a\in M_r(\mathcal V)$ with $L_r(a)\le
1$ such that $\|\ll\varphi, a\gg-\ll\psi, a\gg\|> c-\epsilon$. Let $b=\left[
\begin{array}{cc}
0&a\\
a^*&0
\end{array}\right]$.
We have that $b=b^*\in M_{2r}(\mathcal{V})$, and
\[\begin{array}{rcl}
L_{2r}(b)&=&L_{2r}\left(\left[
  \begin{array}{cc}
    a & 0 \\
    0 & a^*
  \end{array}\right]\left[
  \begin{array}{cc}
    0 & 1 \\
    1 & 0
  \end{array}\right]
\right)\le L_{2r}\left(\left[
  \begin{array}{cc}
    a & 0 \\
    0 & a^*
  \end{array}\right]
\right)\\
&=&\max\{L_r(a), L_r(a^*)\}=L_r(a)\le 1
\end{array}\]
and
\[\begin{array}{rcl}
&&\left\|\ll\varphi, b\gg-\ll\psi, b\gg\right\|\\
&=&\left\|\left[
\begin{array}{cc}
0&\ll\varphi, a\gg-\ll\psi, a\gg\\
(\ll\varphi, a\gg-\ll\psi, a\gg)^\ast&0
\end{array}\right]\right\|\\
&=&\left\|\left[
\begin{array}{cc}
\ll\varphi, a\gg-\ll\psi, a\gg&0\\
0&(\ll\varphi, a\gg-\ll\psi, a\gg)^\ast
\end{array}\right]\right\|\\
&=&\max\{\|\ll\varphi, a\gg-\ll\psi, a\gg\|, \|(\ll\varphi, a\gg-\ll\psi,
a\gg)^\ast\|\}\\
&=&\|\ll\varphi, a\gg-\ll\psi, a\gg\|\\
&>&c-\epsilon.
\end{array}\]
By the arbitrariness of
$\epsilon$,  we obtain that $\sup\{\|\ll \varphi, a\gg -\ll\psi,
a\gg\|:\, a=a^*\in M_r(\mathcal V), L_r(a)\le 1, r\in{\mathbb
N}\}=c$. If $D_{L_n}(\varphi, \psi)=+\infty$, it can be proved
similarly.
\end{proof}

\begin{proposition}\label{pro:34} Let ${\mathcal L}=(L_n)$ be a matrix Lipschitz seminorm on the matrix
order unit space $(\mathcal V, 1)$. Then ${\mathcal L}_{{\mathcal
D}_{\mathcal L}}= (L_{D_{L_n}})$ is a matrix seminorm on
$(\mathcal V, 1)$ and satisfies that
\[L_{D_{L_n}}(a^*)=L_{D_{L_n}}(a)\ \ \hbox{ and }\ \ L_{D_{L_n}}(a)\le L_n(a),\]
for all $n\in\mathbb N$ and $a\in M_n(\mathcal V)$.
\end{proposition}

\begin{proof} For any $v\in M_n(\mathcal V)$, $w\in M_m(\mathcal V)$, $\alpha \in M_{m, n}$
and $\beta \in M_{n, m}$,
we have
\[\begin{array}{rcl}
&&L_{D_{L_{n+m}}}(v\oplus w)\\
&=&\sup\left\{\frac{\|\ll \varphi, v\oplus w\gg -\ll\psi, v\oplus w\gg\|}{D_{L_r}(\varphi,
\psi)}:\,\varphi, \psi\in CS_r(\mathcal V), \varphi\neq \psi, r\in{\mathbb N}\right\}\\
&=&\sup\left\{\frac{\|\ll \varphi, v\gg\oplus\ll\varphi, w\gg -\ll\psi, v\gg\oplus\ll\psi,
w\gg\|}{D_{L_r}(\varphi, \psi)}:\,\varphi, \psi\in CS_r(\mathcal V), \varphi\neq \psi,
r\in{\mathbb N}\right\}\\
&=&\sup\left\{\frac{\|(\ll\varphi, v\gg-\ll\psi, v\gg)\oplus(\ll\varphi, w\gg -\ll\psi,
w\gg)\|}{D_{L_r}(\varphi, \psi)}:\,\varphi, \psi\in CS_r(\mathcal V), \varphi\neq \psi,
r\in{\mathbb N}\right\}\\
&=&\sup\Big\{\frac{\max\{\|\ll\varphi, v\gg-\ll\psi, v\gg\|,
\|\ll\varphi, w\gg -\ll\psi,
w\gg\|\}}{D_{L_r}(\varphi, \psi)}:\\
&&\varphi, \psi\in CS_r(\mathcal V), \varphi\neq \psi, r\in{\mathbb N}\Big\}\\
&=&\max\{L_{D_{L_n}}(v), L_{D_{L_m}}(w)\},
\end{array}\]
and
\[\begin{array}{rcl}
&&L_{D_{L_m}}(\alpha v\beta)=\sup\left\{\frac{\|\ll \varphi,
\alpha v\beta\gg -\ll\psi, \alpha v\beta\gg\|}{D_{L_r} (\varphi,
\psi)}:\,\varphi, \psi\in CS_r(\mathcal V), \varphi\neq \psi,
r\in{\mathbb
N}\right\}\\
&=&\sup\left\{\frac{\|(\alpha\otimes 1_r)(\ll \varphi, v\gg-\ll\psi, v\gg)(\beta\otimes
1_r)\| }{D_{L_r}(\varphi, \psi)}:\,\varphi, \psi\in CS_r(\mathcal V), \varphi\neq \psi,
r\in{\mathbb N}\right\}\\
&\le&\|\alpha\|L_{D_{L_n}}(v)\|\beta\|.
\end{array}\]
So $\mathcal{L_{D_L}}$ is a matrix gauge on $\mathcal V$. For $v\in M_n(\mathcal V)$, we have
\[\begin{array}{rcl}
&&L_{D_{L_n}}(v^*)=\sup\left\{\frac{\|\ll \varphi, v^*\gg
-\ll\psi,  v^*\gg\|}{D_{L_r}(\varphi, \psi)}:
\,\varphi, \psi\in CS_r(\mathcal V), \varphi\neq \psi, r\in{\mathbb N}\right\}\\
&=&\sup\left\{\frac{\|\ll \varphi, v\gg^* -\ll\psi,
v\gg^*\|}{D_{L_r}(\varphi, \psi)}: \,\varphi, \psi\in
CS_r(\mathcal V), \varphi\neq \psi, r\in{\mathbb
N}\right\}=L_{D_{L_n}}(v).
\end{array}\]

If $a=[\lambda_{ij}1]$ for some $[\lambda_{ij}]\in M_n$, then
\[\begin{array}{rcl}
L_{D_{L_n}}(a)&=&\sup\left\{\frac{\|\ll \varphi, [\lambda_{ij}1]\gg -\ll\psi, [\lambda_{ij}1]
\gg\|}{D_{L_r}(\varphi, \psi)}:\,\varphi, \psi\in CS_r(\mathcal V), \varphi\neq \psi,
r\in{\mathbb N}\right\}\\
&=&0=L_n(a),
\end{array}\]
since $\varphi(1)=\psi(1)=1_r$. Suppose $a\in
M_n(\mathcal{V})\setminus M_n({\mathbb C1})$. For $\varphi,
\psi\in CS_r(\mathcal V)$, we have
\[\begin{array}{rcl}
\|\ll\varphi, a\gg-\ll\psi, a\gg\|
&=&L_n(a)\big\|\ll\varphi, \frac{a}{L_n(a)}\gg-\ll\psi, \frac{a}{L_n(a)}\gg\big\|\\
&\le&L_n(a) D_{L_r}(\varphi, \psi).
\end{array}\]
Since this is true for any $\varphi, \psi\in CS_r(\mathcal V)$, we get that $L_{D_{L_n}}
(a)\le L_n(a)$. So $L_{D_{L_n}}(a)\le L_n(a)$ for all $a\in M_n(\mathcal V)$. This also
indicates that $L_{D_{L_n}}(a)<+\infty$ for all $a\in M_n(\mathcal V)$. Therefore,
$\mathcal{L_{D_L}}$ is a matrix seminorm on  $\mathcal V$.
\end{proof}

\begin{corollary}\label{co:35} Let ${\mathcal L}=(L_n)$ be a matrix Lipschitz seminorm on the matrix order
unit space $(\mathcal V, 1)$. Then the graded set
\[\mathcal{L}^{1}_{\mathcal{D_L}}=(L_{D_{L_n}}^1),\]
where $L_{D_{L_n}}^1=\{a\in M_n(\mathcal V):\,L_{D_{L_n}}(a)\le
1\}$, is absolutely matrix convex.
\end{corollary}

Let $V$ be a locally convex vector space, and let $L$ be a gauge on $V$. Then $L$ is {\it
lower semicontinuous} if for any net $\{a_\lambda\}$ in $V$ which converges to $a\in V$ we
have $L(a)\le\liminf_{\lambda}\{L(a_\lambda)\}$. Equivalently, for one, hence every,
$t\in\mathbb R$ with $t>0$, the set
\[\{a\in V:\,L(a)\le t\}\]
is closed in $V$. If $V$ and $W$ are two locally convex vector
spaces in duality, we have that $M_n(V)$ and $M_m(W)$ are in
duality under the matrix pairing. The {\it weak topology} on
$M_n(V)$ is defined to be that determined by the matrix pairing.

\begin{proposition}\label{pro:36} Let $(\mathcal V, 1)$ be a matrix order unit space, and let ${\mathcal
L}=(L_n)$ be any matrix Lipschitz seminorm on $(\mathcal V, 1)$. Then ${\mathcal L_{\mathcal
D_\mathcal L}}=(L_{D_{L_n}})$ is lower semicontinuous, that is, each $L_{D_{L_n}}$  is lower
semicontinuous.
\end{proposition}

\begin{proof} Suppose $\{a_i\}$ is a net in $M_n(\mathcal V)$ which converges weakly to
$a\in M_n(\mathcal V)$. For $r\in\mathbb N$, $\varphi, \psi\in CS_r(\mathcal V)$ with
$\varphi \neq\psi$,  we have
\[\left\|\frac{\ll\varphi , a\gg-\ll\psi , a\gg}{D_r(\varphi,  \psi)}-\frac{\ll\varphi ,
a_i\gg-\ll\psi , a_i\gg}{D_r(\varphi,
\psi)}\right\|=\frac{\|\ll\varphi -\psi , a-a_i\gg\|}{D_r(\varphi,
\psi)}\] is small when $i$ is large, and therefore, for each
$\epsilon >0$,
\[\bigg\|\frac{\ll\varphi , a\gg-\ll\psi , a\gg}{D_r(\varphi,  \psi)}\bigg\|\le\bigg
\|\frac{\ll\varphi , a_i\gg-\ll\psi , a_i\gg}{D_r(\varphi,  \psi)}\bigg\|+\epsilon,\]
for $i$ sufficiently large. It follows that
\[\begin{array}{rcl}
&&\big\|\frac{\ll\varphi , a\gg-\ll\psi , a\gg}{D_r(\varphi,  \psi)}\big\|\\
&\le&\sup\left\{\frac{\|\ll\varphi , a_i\gg-\ll\psi , a_i\gg\|}{D_r(\varphi,  \psi)}:\,
\varphi,  \psi\in CS_r(\mathcal V), \varphi\neq\psi, r\in{\mathbb N}\right\}+\epsilon,
\end{array}\]
for $i$ sufficiently large, and hence that
\[\frac{\|\ll\varphi , a\gg-\ll\psi , a\gg\|}{D_r(\varphi,  \psi)}\le\liminf_i L_{D_{L_n}}
(a_i)+\epsilon.\] Let $\epsilon$ tend to $0$, and then we get that
$L_{D_{L_n}}(a)\le\liminf_i L_{D_{L_n}}(a_i)$. Therefore,
$L_{D_{L_n}}$  is lower  semicontinuous.
\end{proof}

\section{Recovering $\mathcal L$ from $\mathcal D_{\mathcal L}$}\label{sec:4}

Let $(\mathcal V, 1)$ be a matrix order unit space and let
${\mathcal L}=(L_n)$ be a matrix Lipschitz seminorm on $(\mathcal
V, 1)$. We denote by $\|\cdot\|=(\|\cdot\|_n)$ the matrix norm
determined by the matrix order on $(\mathcal V, 1)$. Since
${\mathbb C}1$ is a closed subspace of $\mathcal V$, $M_n({\mathbb
C}1)$ is closed in $M_n(\mathcal V)$ for each $n\in\mathbb
N$ by Proposition 2.1(2) in \cite{ru}. We may therefore use the identification
\[M_n(\mathcal V/({\mathbb C}1))\cong M_n(\mathcal V)/M_n({\mathbb C}1)\]
to define a matrix norm $\|\cdot\|^{\sim}=(\|\cdot\|_n^{\sim})$ on
the quotient space $\tilde{\mathcal V}=\mathcal V/({\mathbb C}1)$
from $\|\cdot\|= (\|\cdot\|_n)$ by Theorem 4.2 in \cite{ru}(or see \cite{efruan}). But in addition
to this matrix norm, by identifying $M_n(\mathcal V/({\mathbb
C}1))$ with $M_n(\mathcal V)/M_n({\mathbb C}1)$, we may also
obtain a sequence of quotient seminorm $\widetilde{\mathcal
L}=(\widetilde{L}_n)$ from $\mathcal L$. $\widetilde{\mathcal L}$
is also a matrix norm on $\tilde{\mathcal V}$ since the null space
of each $L_n$ is $M_n({\mathbb C}1)$(see page 173 in \cite{efwe}). It is clear that
$\tilde{\mathcal V}$ is also self-adjoint since ${\mathbb C}1$ is
self-adjoint.

From Proposition 5.2 in \cite{wu}, we have

\begin{lemma}\label{le:41} Let $(\mathcal V, 1)$ be a matrix order unit space. For arbitrary $r\in\mathbb N$, we denote
\[CB_r^2(\tilde{\mathcal{V}})=\{f\in M_r((\tilde{\mathcal V})^*):\, \|f\|_{cb}\le 2\}\]
and
\[SB_r^2(\tilde{\mathcal{V}})=\{f\in M_r((\tilde{\mathcal V})^*):
\, f^*=f,\,\|f\|_{cb}\le 2\}.\]
Then
\[CB_r^2(\tilde{\mathcal{V}})\subseteq\{\varphi_1-\varphi_2+i(\varphi_3-\varphi_4):
\, \varphi_i\in CS_r(\mathcal V),
i=1, 2, 3, 4\}\]
and
\[SB_r^2(\tilde{\mathcal{V}})=\{\varphi_1-\varphi_2:\, \varphi_i\in CS_r(\mathcal V),
i=1, 2\}.\] Here we view $M_r((\tilde{\mathcal{V}})^*)$ as the
subspace of $M_r(\mathcal{V}^*)$ consisting of those $f\in
M_r(\mathcal{V}^*)$ such that $f(a)=0_r$ for $a\in\mathbb{C}1$.
\end{lemma}

\begin{lemma}\label{le:42} Let $(\mathcal V, 1)$ be a matrix order unit space. Then for any $f\in
M_r((\tilde{\mathcal{V}})^{*})$, we have
\[\|f^{*}\|_{cb}=\|f\|_{cb}.\]
\end{lemma}

\begin{proof} Let $\|\cdot\|=(\|\cdot\|_n)$ be the matrix
norm on $(\mathcal V, 1)$ determined by the matrix order on it. By Proposition 2.1 in \cite{ble}, 
we have
\[\begin{array}{rcl}
\|f^*\|_{cb}&=&\sup\{\|\ll f^*, a\gg\|: a=[a_{ij}]\in
M_r(\mathcal{V}), \|a\|_r\le 1\}\\
&=&\sup\{\|[f^*(a_{ij})]\|: a=[a_{ij}]\in
M_r(\mathcal{V}), \|a\|_r\le 1\}\\
&=&\sup\{\|[f(a_{ij}^\ast)^\ast]\|: a=[a_{ij}]\in
M_r(\mathcal{V}), \|a\|_r\le 1\}\\
&=&\sup\{\|[f(a_{ji}^\ast)]\|: a=[a_{ij}]\in
M_r(\mathcal{V}), \|a\|_r\le 1\}\\
&=&\sup\{\|\ll f, [a_{ji}^\ast]\gg\|: a=[a_{ij}]\in
M_r(\mathcal{V}), \|a\|_r\le 1\}\\
&=&\sup\{\|\ll f, [a_{ij}]\gg\|: a=[a_{ij}]\in
M_r(\mathcal{V}), \|a\|_r\le 1\}\\
&=&\|f\|_{cb}.
\end{array}\]
\end{proof}

For a matrix Lipschitz seminorm ${\mathcal L}=(L_n)$ on a matrix
order unit space $(\mathcal V, 1)$ and $t>0$, denote by
${\mathcal{L}^t}$ the sequence of sets
\[L^t_n=\{a\in M_n(\mathcal V):\,L_n(a)\le t\},\]
and by $\overline{\mathcal{L}}^t=(\overline{L}^t_n)$ the norm
closure of ${\mathcal{L}^t}$, that is, each $\overline{L}^t_n$ is
the norm closure of $L^t_n$.

\begin{lemma}\label{le:43}
Let $(\mathcal V, 1)$ be a matrix order unit space, and let
${\mathcal L}= (L_n)$ be a matrix Lipschitz seminorm on $(\mathcal
V, 1)$. For $n\in\mathbb{N}$, define $L_n^{'}:
M_n(\mathcal{V}^*)\longmapsto [0,+\infty]$ by
\[L^{'}_n(f)=\sup\{\|\ll f, a\gg\|:\,a\in L_r^1, r\in{\mathbb
N}\},\ \ \ f\in M_n(\mathcal{V}^*).\] Then ${\mathcal
L}^{'}=(L^{'}_n)$ is a matrix gauge on ${\mathcal V}^*$, and
\[L_n^{'}(f)=\sup\{\|\ll f, a\gg\|:\, a\in L_r^1, a=a^*,
r\in\mathbb N\},\]for $f\in M_n(\mathcal{V}^*)$, and
\[L^{'}_n(\varphi-\psi)=D_{L_n}(\varphi, \psi),\] for $\varphi,
\psi\in CS_{n}(\mathcal{V})$.
\end{lemma}

\begin{proof} It is easy to verify that ${\mathcal
L}^{'}=(L^{'}_n)$ is a matrix gauge on $\mathcal{V}^\ast$.

For arbitrary $r\in \mathbb N$ and $a\in L_r^1$, we have that
$L_r(a^*)=L_r(a)\le 1$. Since
\[L_{2r}\bigg(\left[\begin{array}{cc}
0&a\\
a^*&0\end{array}\right]\bigg)= L_{2r}\bigg(\left[\begin{array}{cc}
a&0\\
0&a^*\end{array}\right]\left[\begin{array}{cc}
0&1\\
1&0\end{array}\right]\bigg)\le L_{2r}\bigg(\left[\begin{array}{cc}
a&0\\ 0&a^*\end{array} \right]\bigg),\] we get that
$L_{2r}\bigg(\left[\begin{array}{cc}
0&a\\
a^*&0\end{array}\right]\bigg)\le
L_{2r}\bigg(\left[\begin{array}{cc}
a&0\\
0&a^*\end{array}\right]\bigg)$. Similarly, we have
\[L_{2r}\bigg(\left[\begin{array}{cc}
a&0\\
0&a^*\end{array}\right]\bigg)\le
L_{2r}\bigg(\left[\begin{array}{cc}
0&a\\
a^*&0\end{array}\right]\bigg).\] So
\[L_{2r}\bigg(\left[\begin{array}{cc}
0&a\\
a^*&0\end{array}\right]\bigg)= L_{2r}\bigg(\left[\begin{array}{cc}
a&0\\
0&a^*\end{array}\right]\bigg)=L_r(a)\le 1.\] Given $f\in
M_n(\mathcal{V}^*)$. We have
\[\begin{array}{rcl}
&&\left\|\ll f, \left[\begin{array}{cc}
0&a\\
a^*&0\end{array}\right]\gg\right\|=\left\|\left[\begin{array}{cc}
0&\ll f, a\gg\\
\ll f, a^*\gg&0\end{array}\right]\right\|\\
&=&\left\|\left[\begin{array}{cc}
\ll f, a\gg&0\\
0&\ll f, a^*\gg\end{array}\right]\right\|\ge\|\ll f,
a\gg\|.\end{array}\] Then $L_n^{'}(f)\le\sup\{\|\ll f,
a\gg\|:\,a\in L_r^1, a^*=a, r\in\mathbb N\}$. Therefore,
\[L_n^{'}(f)=\sup\{\|\ll f, a\gg\|:\, a\in L_r^1, a=a^*,
r\in\mathbb N\}.\] For $\varphi, \psi\in CS_n(\mathcal{V})$, we
have
\[\begin{array}{rcl}
L^{'}_n(\varphi-\psi)
&=&\sup\{\|\ll\varphi-\psi, a\gg\|:\,a\in L_r^1, r\in{\mathbb N}\}\\
&=&\sup\{\|\ll\varphi, a\gg-\ll\psi, a\gg\|:\,a\in L_r^1,
r\in{\mathbb
N}\}\\
&=&D_{L_n}(\varphi, \psi).
\end{array}\]
\end{proof}

\begin{theorem}\label{th:44} Let $(\mathcal V, 1)$ be a matrix order unit space, and let ${\mathcal L}=
(L_n)$ be a lower semicontinuous matrix Lipschitz seminorm on $(\mathcal V, 1)$. Then
\[\mathcal L_{\mathcal D_\mathcal L}=\mathcal L.\]
\end{theorem}

\begin{proof}
The generalized bipolar theorem says that $({\mathcal
L}^1)^{\circledcirc\circledcirc}$ is the smallest weakly closed
absolutely matrix convex set containing $\mathcal L^1$(Proposition 4.1 in \cite{efwe}). 
Because the closure of a convex set in a locally convex topological vector space coincide with the closure under the 
weak topology, the bipolar theorem says also that the absolute operator bipolar of $\mathcal{L}^1$ is the smallest norm closed 
absolutely matrix convex set containing $\mathcal{L}^1$.
Since ${\mathcal L}=(L_n)$ is a matrix gauge, ${\mathcal L}^1$ is
absolutely matrix convex. The lower semicontinuity of $\mathcal L$
implies that $\mathcal L^1$ is norm closed. Thus
\[({\mathcal L}^1)^{\circledcirc\circledcirc}=\mathcal L^1.\]
But for $n\in\mathbb N$, we have
\[\begin{array}{rcl}
&&(L_n^1)^{\circledcirc}\\
&=&\{f\in M_n(\mathcal V^*):\,\|\ll f, a\gg\|\le 1\hbox{ for all
}a\in M_r(\mathcal V),
L_r(a)\le 1, r\in{\mathbb N}\}\\
&=&\{f\in M_n(\mathcal V^*):\, L^{'}_n(f)\le 1\}
\end{array}\]
and
\[\begin{array}{rcl}
&&(L_n^1)^{\circledcirc\circledcirc}\\
&=&\{f\in M_n(\mathcal V^*):\, L^{'}_n(f)\le 1\}^\circledcirc\\
&=&\{a\in M_n(\mathcal V):\,\|\ll f, a\gg\|\le 1\hbox{ for all
}f\in M_r(\mathcal V^*), L^{'}_r(f)\le 1, r\in{\mathbb N}\}.
\end{array}\]

Suppose $a\in M_n(\mathcal V)$ with $L_{D_{L_n}}(a)\le 1$. Then
$\|\ll\varphi, a\gg-\ll\psi, a\gg\|\le D_{L_r}(\varphi, \psi)$ for
all $\varphi, \psi\in CS_r(\mathcal V)$ and arbitrary $r\in\mathbb
N$. So by Lemma \ref{le:43}, we have that $\|\ll\varphi,
a\gg-\ll\psi, a\gg\|\le L_r^{'}(\varphi- \psi)$ for all $\varphi,
\psi\in CS_r(\mathcal V)$ and arbitrary $r\in\mathbb N$. Thus for
$f\in SB^2_r(\tilde{\mathcal{V}})$ and $r\in\mathbb{N}$, we have
\[\|\ll f, a\gg\|\le L_r^{'}(f)\]
by Lemma \ref{le:41}.

Given $f\in CB^2_r(\tilde{\mathcal{V}})$. We have
\[\begin{array}{rcl}
\bigg\|\left[\begin{array}{cc}
0&f\\
f^*&0\end{array}\right]\bigg\|_{cb}
&=&\bigg\|\left[\begin{array}{cc}
f&0\\
0&f^\ast\end{array}\right]
\left[\begin{array}{cc}
0&1\\
1&0\end{array}\right]\bigg\|_{cb}\\
&\le&\bigg\|\left[\begin{array}{cc}
f&0\\
0&f^\ast\end{array}\right]\bigg\|_{cb}\\
&=&\max\{\|f\|_{cb},
\|f^\ast\|_{cb}\}\\
&=&\|f\|_{cb}\le 2\end{array}\] 
by Lemma \ref{le:42}. Thus
\[\begin{array}{rcl}
\|\ll f, a\gg\|&\le&\bigg\|\ll\left[\begin{array}{cc}
0&f\\
f^*&0\end{array}\right], a\gg\bigg\|\\
&\le&L_{2r}^{'}\bigg(\left[\begin{array}{cc}
0&f\\
f^*&0\end{array}\right]\bigg)=L_r^{'}(f),\end{array}\] 
that is, if $f\in M_r(\mathcal V^*)$ with $L^{'}_r(f)\le 1$, then
$\|\ll f, a\gg\|\le 1$. Now, we have that $a\in
(L_n^1)^{\circledcirc\circledcirc}$. So $a\in L_n^1$. Therefore,
$L_n(a)\le L_{D_{L_n}} (a)$ for all $a\in M_n(\mathcal V)$.
Combining this and Proposition \ref{pro:34}, we obtain that $L_n=
L_{D_{L_n}}$. So $\mathcal{L_{D_L}=L}$.
\end{proof}

\begin{corollary}\label{co:45}
Let $(\mathcal{V},1)$ be a matrix order unit space, and let
$\mathcal{L}=(L_n)$ be a matrix Lipschitz seminorm on
$(\mathcal{V},1)$. Then $\mathcal{L_{D_L}}$ is the largest lower
semicontinuous matrix Lipschitz seminorm smaller than
$\mathcal{L}$, and $\mathcal{D_{L_{D_L}}}=\mathcal{D_L}$.
\end{corollary}

\section{Matrix norms from matrix metrics}\label{sec:5}

In this section, we discuss the relation of the matrix metrics on
$\mathcal{CS(V)}$ and the matrix norms on
$(\tilde{\mathcal{V}})^*$. Firstly, we make:

\begin{definition}\label{def:51} Let $V$ be a vector space and let $\mathbf K=(K_n)$ be a graded set with
$K_n\subseteq M_n(V)$. A {\it matrix metric} ${\mathcal D}=(D_n)$ on $\mathbf K$ is a sequence
of metrics
\[D_n:\, K_n\times K_n\longmapsto [0, +\infty)\]
such that
\begin{enumerate}
\item if $x, u\in K_m$ and $y, v\in K_n$ such that $x\oplus y, u\oplus v\in K_{m+n}$, then
$D_{m+n}(x\oplus y, u\oplus v)=\max\{D_m(x, u), D_n(y, v)\}$;
\item if $x, u\in K_m$ and $\alpha\in M_{m, n}$ with $\alpha^*\alpha=1_n$ such that
$\alpha^* x\alpha, \alpha^* u\alpha\in K_n$, then $D_n(\alpha^*
x\alpha, \alpha^* u\alpha) \le D_m(x, u)$.
\end{enumerate}
\end{definition}

\begin{example}\label{ex:52}
Assume $\mathcal L=(L_n)$ is a matrix Lipschitz seminorm on the
matrix order unit space $(\mathcal V, 1)$ and the image of $L_1^1$ in $\tilde{\mathcal{V}}$ is totally
bounded for $\|\cdot\|_1^\sim$. Then the sequence of
metrics $D_{L_n}$ defined in section \ref{sec:3} is a matrix
metric on $\mathcal{CS(V)}$(see Theorem 5.3 in \cite{wu}).
\end{example}

Suppose $K$ is a convex subset of a vector space $V$. A metric $\rho$ on $K$ is {\it convex}
if for every $x, y, z\in K$ and $t\in [0, 1]$ we have
\[\rho(x, ty+(1-t)z)\le t\rho(x, y)+(1-t)\rho(x, z).\]
A metric $\rho$ on $K$ is {\it midpoint balanced} if whenever $x, y, u, v\in K$ and
\[\frac{x+v}{2}=\frac{y+u}{2}\]
holds, it follows that $\rho(x, u)=\rho(y, v)$. A metric $\rho$ on
$K$ is {\it midpoint concave} if for any $x, y, u, v\in K$ we have
\[\rho\Big(\frac{x+u}{2}, \frac{y+v}{2}\Big)\le\frac12(\rho(x, y)+\rho(u, v))\]
\cite{ri}.

A {\it matrix convex set} in a vector space $V$ is a collection
$\mathbf{K}=(K_n)$ of subsets $K_n\subseteq M_n(V)$ such that
\[\sum_{i=1}^{k}\gamma_i^*v_i\gamma_i\in K_n\]
for all $v_i\in K_{n_i}$ and $\gamma_i\in M_{n_i,n}$ for
$i=1,2,\cdots,k$ satisfying
$\sum_{i=1}^k\gamma_i^*\gamma_i=1_n$\cite{wewi, witt, efwi}.

\begin{definition}\label{def:53} We say that a matrix metric $\mathcal D=(D_n)$ on a matrix convex set
$\mathbf K$ in a vector space $V$ is {\it convex} if each $D_n$ is
convex; $\mathcal D$ is {\it midpoint balanced} if each $D_n$ is
midpoint balanced; $\mathcal D$ is {\it midpoint concave} if each
$D_n$ is midpoint concave.
\end{definition}

\begin{lemma}\label{le:54}
Let $(\mathcal V, 1)$ be a matrix order unit space, and let
${\mathcal D}= (D_n)$ be a matrix metric on  ${\mathcal C\mathcal
S}(\mathcal V)$. If $\mathcal D$ is convex, midpoint balanced and
midpoint concave, then there is a norm $Q_n^{''}$ on each real
vector space $\mathbb{R}SB_n^2(\tilde{\mathcal{V}})$ such that
\[Q_n^{''}(\alpha^*f\alpha)\le\|\alpha\|^2Q_m^{''}(f),\]
\[Q_{m+n}^{''}(f\oplus g)=\max\{Q_m^{''}(f),Q_n^{''}(g)\},\]
for $f\in\mathbb{R}SB_m^2(\tilde{\mathcal{V}})$,
$g\in\mathbb{R}SB_n^2(\tilde{\mathcal{V}})$ and $\alpha\in
M_{m,n}$, and
\[Q_n^{''}(\varphi-\psi)=D_n(\varphi, \psi),\]
for all $\varphi, \psi\in CS_n(\mathcal{V})$.
\end{lemma}

\begin{proof} By Lemma \ref{le:41}, we have that
$SB_n^2(\tilde{\mathcal{V}})=CS_n(\mathcal V)-CS_n(\mathcal
V)(n\in\mathbb{N})$. For $f=\varphi-\psi\in
SB_n^2(\tilde{\mathcal{V}})$, where $\varphi, \psi\in
CS_n(\mathcal V)$, set
\[Q_n^{'}(f)=D_n(\varphi, \psi).\]
Suppose $\varphi_1, \varphi_2, \psi_1, \psi_2\in CS_n(\mathcal V)$
and $\varphi_1-\psi_1= \varphi_2-\psi_2$. Then
$\frac{\varphi_1+\psi_2}{2}=\frac{\varphi_2+\psi_1}{2}$, and hence
$D_n(\varphi_1, \psi_1)=D_n(\varphi_2, \psi_2)$ since $D_n$ is
midpoint balanced. So $Q_n^{'}$ is well-defined.

If $t\in [-1, 1]$ and $f=\varphi-\psi\in
SB_n^2(\tilde{\mathcal{V}})$ with $\varphi, \psi\in CS_n(\mathcal
V)$, we have that $tf=\mathrm{sgn}(t)\abs{t}(\varphi-\psi)=\mathrm{sgn}(t)[\varphi-(\abs{t}\psi+(1-\abs{t}) \varphi)]$. So $tf\in
SB_n^2(\tilde{\mathcal{V}})$ by the convexity of $CS_n(\mathcal
V)$. For $t\in [0, 1]$, since
\[\begin{array}{rcl}
D_n(\varphi, \psi)&\le& D_n(\varphi, t\psi+(1-t)\varphi)+D_n(t\psi+(1-t)\varphi, \psi)\\
&\le& t D_n(\varphi, \psi)+(1-t)D_n(\varphi, \psi)=D_n(\varphi,
\psi)\end{array}\] by the convexity of $D_n$, it follows that
$tD_n(\varphi, \psi)=D_n(\varphi, t\psi+(1-t)\varphi)$. Therefore,
$Q_n^{'}(tf)=D_n(\varphi, t\psi+(1-t)\varphi)=tD_n(\varphi,
\psi)=tQ_n^{'}(f)$, and $Q_n^{'}(-tf)=Q_n^{'}(tf)=tQ_n^{'}(f)$. So
\[Q_n^{'}(tf)=\abs{t}Q_n^{'}(f),\, \abs{t}\le 1.\]
Now for $f=\varphi_1-\psi_1, g=\varphi_2-\psi_2\in
SB_n^2(\tilde{\mathcal{V}})$ with $\varphi_1, \varphi_2, \psi_1,
\psi_2\in CS_n(\mathcal V)$ and $f+g\in
SB_n^2(\tilde{\mathcal{V}})$, we obtain
\[\begin{array}{rcl}
Q_n^{'}(f+g)&=&2Q_n^{'}\big(\frac{f+g}{2}\big)=
2Q_n^{'}\big(\frac{\varphi_1+\varphi_2-\psi_1-\psi_2}{2}\big)\\
&=&2D_n\big(\frac{\varphi_1+\varphi_2}{2},
\frac{\psi_1+\psi_2}{2}\big)
\le D_n(\varphi_1, \psi_1)+D_n(\varphi_2, \psi_2)\\
&=&Q_n^{'}(f)+Q_n^{'}(g)\end{array}\] by the midpoint concavity of
$D_n$. For $f_1=\varphi_1-\psi_1\in SB_m^2(\tilde{\mathcal{V}})$
with $\varphi_1, \psi_1\in CS_m(\mathcal V)$ and
$f_2=\varphi_2-\psi_2\in SB_n^2(\tilde{\mathcal{V}})$ with
$\varphi_2, \psi_2\in CS_n(\mathcal V)$, we have
\[\begin{array}{rcl}
Q_{m+n}^{'}(f_1\oplus f_2)&=&Q_{m+n}^{'}(\varphi_1\oplus\varphi_2-\psi_1\oplus\psi_2)\\
&=&D_{m+n}(\varphi_1\oplus\varphi_2, \psi_1\oplus\psi_2)\\
&=&\max\{D_m(\varphi_1, \psi_1), D_n(\varphi_2, \psi_2)\}\\
&=&\max\{Q_m^{'}(f_1), Q_n^{'}(f_2)\}.\end{array}\]

Given $f=\varphi-\psi\in SB_m^2(\tilde{\mathcal{V}})$ with
$\varphi, \psi\in CS_m(\mathcal V)$ and $\alpha\in M_{m, n}$ with
$\alpha^*\alpha\le 1_n$. Choose $\phi\in CS_n(\mathcal{V})$. Since
\[
\left[\begin{array}{cc} f&0\\
0&0_n\end{array}\right]=\left[\begin{array}{cc}
\varphi-\psi&0\\
0&0_n\end{array}\right]=\left[\begin{array}{cc}
\varphi&0\\
0&\phi\end{array}\right]-\left[\begin{array}{cc}
\psi&0\\
0&\phi\end{array}\right]=\varphi_1-\psi_1,\] where
$\varphi_1=\left[\begin{array}{cc}
\varphi&0\\
0&\phi\end{array}\right]$ and $\psi_1=\left[\begin{array}{cc}
\psi&0\\
0&\phi\end{array}\right]$, and $\varphi_1, \varphi_2\in
CS_{m+n}(\mathcal{V})$, we see that $\left[\begin{array}{cc}
f&0\\
0&0_n\end{array}\right]\in SB_{m+n}^2(\tilde{\mathcal{V}})$. So
\[\alpha^*f\alpha=\left[\begin{array}{c}
\alpha\\
\sqrt{1_n-\alpha^*\alpha}\end{array}\right]^*\left[\begin{array}{cc}
f&0\\
0&0_n\end{array}\right]\left[\begin{array}{c}
\alpha\\
\sqrt{1_n-\alpha^*\alpha}\end{array}\right]\in
SB_n^2(\tilde{\mathcal{V}})\] by the matrix convexity of
$\mathcal{CS(V)}$, and hence
\[\begin{array}{rcl}
Q_n^{'}(\alpha^*f\alpha)
&=&Q_n^{'}\bigg(\beta^*\left[\begin{array}{cc}
f&0\\
0&0_n\end{array}\right]\beta\bigg)=Q_n^{'}\bigg(\beta^*\varphi_1\beta-\beta^*\psi_1\beta\bigg)\\
&=&D_n\bigg(\beta^*\varphi_1\beta, \beta^*\psi_1\beta\bigg)\le
D_{m+n}(\varphi_1, \psi_1)\\
&=&Q_{m+n}^{'}(\varphi_1-\psi_1)=Q_{m+n}^{'}\bigg(\left[\begin{array}{cc}
f&0\\
0&0_n\end{array}\right]\bigg)=Q_m^{'}(f),\end{array}\] where
$\beta=\left[\begin{array}{c}
\alpha\\
\sqrt{1_n-\alpha^*\alpha}\end{array}\right]$.

Denote
\[W_n=\mathbb{R}SB_n^2(\tilde{\mathcal{V}}), \, n\in\mathbb N.\]
For $f\in W_n$, there is a $t>0$ such that $tf\in
SB_n^2(\tilde{\mathcal{V}})$. Set
\[Q_n^{''}(f)=t^{-1}Q_n^{'}(tf).\]
Then it is clear that $Q_n^{''}$ is well-defined and
\[Q_n^{''}(sf)=\abs{s}Q_n^{''}(f),\]
for $s\in\mathbb{R}$ and $f\in W_n$. If
$f=t_1(\varphi_1-\psi_1)\in W_n$ and $g=t_2(\varphi_2-\psi_2)\in
W_n$ with $\varphi_1, \varphi_2, \psi_1, \psi_2\in CS_n(\mathcal
V)$ and $t_1, t_2>0$, we may assume that $t_1, t_2> 1$, and then
we have
\[\begin{array}{rcl}
Q_n^{''}(f+g)&=&Q_n^{''}(t_1(\varphi_1-\psi_1)+t_2(\varphi_2-\psi_2))\\
&=&Q_n^{''}((t_1\varphi_1+t_2\varphi_2)-(t_1\psi_1+t_2\psi_2))\\
&=&Q_n^{''}\left((t_1+t_2)\left(\frac{t_1\varphi_1+t_2\varphi_2}{t_1+t_2}-\frac{t_1\psi_1+t_2
\psi_2}{t_1+t_2}\right)\right)\\
&=&(t_1+t_2)Q_n^{'}\left(\frac1{t_1+t_2}f+\frac1{t_1+t_2}g\right)\\
&\le&(t_1+t_2)Q_n^{'}\left(\frac1{t_1+t_2}f\right)+(t_1+t_2)Q_n^{'}\left(\frac1{t_1+t_2}g
\right)\\
&=&Q_n^{''}(f)+Q_n^{''}(g),
\end{array}\]
and for $\alpha\in M_{n, m}$, we have
\[\begin{array}{rcl}
&&Q_m^{''}(\alpha^*f\alpha)=Q_m^{''}\left(t_1\|\alpha\|^2\left(\frac{\alpha}
{\|\alpha\|}\right)^*
\frac{f}{t_1}\left(\frac{\alpha}{\|\alpha\|}\right)\right)\\
&=&t_1\|\alpha\|^2Q_m^{'}\left(\left(\frac{\alpha}{\|\alpha\|}\right)^*\frac{f}{t_1}
\left(\frac{\alpha}{\|\alpha\|}\right)\right)
\le t_1\|\alpha\|^2Q_n^{'}\left(\frac{f}{t_1}\right)\\
&=&t_1\|\alpha\|^2Q_n^{''}\left(\frac{f}{t_1}\right)=\|\alpha\|^2Q_n^{''}(f)\end{array}\]
if $\alpha\neq 0_{n, m}$ and clearly it holds for $\alpha=0_{n,
m}$. Let $f=t_1f_1\in W_m$ with $t_1> 1$ and $f_1\in
SB_m^2(\tilde{\mathcal{V}})$ and $g=t_2g_1\in W_n$ with $t_2> 1$
and $g_1\in SB_n^2(\tilde{\mathcal{V}})$. Setting $t=\max\{t_1,
t_2\}$, we have
\[\begin{array}{rcl}
&&Q_{m+n}^{''}(f\oplus
g)=Q_{m+n}^{''}\left(t\left(\frac{t_1}{t}f_1\oplus \frac{t_2}tg_1
\right)\right)\\
&=&tQ_{m+n}^{'}\left(\frac{t_1}tf_1\oplus \frac{t_2}tg_1\right)
=t\max\left\{Q_m^{'}\left(\frac{f}t\right), Q_n^{'}\Big(\frac{g}t\Big)\right\}\\
&=&t\max\left\{Q_m^{''}\left(\frac{f}t\right),
Q_n^{''}\Big(\frac{g}t\Big)\right\} =\max\left\{Q_m^{''}(f),
Q_n^{''}(g)\right\}.
\end{array}\]
\end{proof}

\begin{theorem}\label{th:55}
Let $(\mathcal V, 1)$ be a matrix order unit space, and let
${\mathcal D}= (D_n)$ be a matrix metric on  ${\mathcal C\mathcal
S}(\mathcal V)$. Then there is a matrix norm $\mathcal{Q}=(Q_n)$
on $(\tilde{\mathcal{V}})^*$ such that
\[Q_n(\varphi-\psi)=D_n(\varphi,\psi),\]
for all $n\in\mathbb{N}$ and $\varphi,\psi\in CS_n(\mathcal{V})$,
and \[Q_n(f^*)=Q_n(f),\] for each $f\in
M_n((\tilde{\mathcal{V}})^*)$ and all $n\in\mathbb{N}$, if and
only if $\mathcal D$ is convex, midpoint balanced, and midpoint
concave. The matrix norm $\mathcal{Q}$
is unique.
\end{theorem}

\begin{proof}
The uniqueness follows immediately from the equality
\[Q_n(f)=\max\left\{Q_n(f^*), Q_n(f)\right\}=Q_{2n}\left(\left[\begin{array}{cc}
f^*&0\\
0&f\end{array}\right]\right)=Q_{2n}\left(\left[\begin{array}{cc}
0&f\\
f^*&0\end{array}\right]\right),\] where $f\in
M_n((\tilde{\mathcal{V}})^*)$.

If $Q_n(\varphi-\psi)=D_n(\varphi,\psi)$, for all $n\in\mathbb{N}$
and $\varphi,\psi\in CS_n(\mathcal{V})$, for some matrix norm
$\mathcal{Q}=(Q_n)$ on $(\tilde{\mathcal{V}})^*$ satisfying
$Q_n(f^*)=Q_n(f)$ for each $f\in M_n((\tilde{\mathcal{V}})^*)$ and
all $n\in\mathbb{N}$, it is easily seen that $\mathcal D$ is
convex, midpoint balanced, and midpoint concave.

Now suppose $\mathcal D$ is convex, midpoint balanced, and midpoint
concave. By Lemma \ref{le:54}, there
is a norm $Q_n^{''}$ on each real vector space
$\mathbb{R}SB_n^2(\tilde{\mathcal{V}})$ such that
$Q_n^{''}(\alpha^*f\alpha)\le\|\alpha\|^2Q_n^{''}(f)$ and
$Q_{m+n}^{''}(f\oplus g)=\max\{Q_m^{''}(f),Q_n^{''}(g)\}$ for
$f\in\mathbb{R}SB_m^2(\tilde{\mathcal{V}})$ and
$g\in\mathbb{R}SB_n^2(\tilde{\mathcal{V}})$ and $\alpha\in
M_{m,n}$, and $Q_n^{''}(\varphi-\psi)=D_n(\varphi, \psi)$ for all
$\varphi, \psi\in CS_n(\mathcal{V})$. Denote
\[U_n=\mathbb{R}SB_n^2(\tilde{\mathcal{V}})+i\mathbb{R}SB_n^2(\tilde{\mathcal{V}}),\,\,
n\in\mathbb N.\] 
For $f\in U_n$, set
\[Q_n(f)=Q_{2n}^{''}\bigg(\left[\begin{array}{cc}
0&f\\
f^*&0\\
\end{array}\right]\bigg).\]
For $f=t(\varphi-\psi)\in\mathbb{R}
SB_n^2(\tilde{\mathcal{V}})$ with $t> 1$ and $\varphi, \psi\in
CS_n(\mathcal V)$, we have that
$\bigg\|\frac{1}{t}\left[\begin{array}{cc}
0&f\\
f&0\\
\end{array}\right]\bigg\|_{cb}=\frac{1}{t}\|f\|_{cb}\le 2$,
that is, $\frac{1}{t}\left[\begin{array}{cc}
0&f\\
f&0\\
\end{array}\right]\in SB^2_{2n}(\tilde{\mathcal{V}})$. So
\[Q_n(f)=Q_{2n}^{''}\bigg(\left[\begin{array}{cc}
0&f\\
f&0\\
\end{array}\right]\bigg)=tQ_{2n}^{'}\bigg(\frac{1}{t}\left[\begin{array}{cc}
0&f\\
f&0\\
\end{array}\right]\bigg).\]
Because ${\mathcal D}= (D_n)$ is a matrix metric and
\[\begin{array}{rcl}
\frac1t\left[\begin{array}{cc}
0&f\\
f&0\\
\end{array}\right]&=&\frac1t\alpha^\ast
\left[\begin{array}{cc}
f&0\\
0&-f\\
\end{array}\right]\alpha\\
&=&\alpha^\ast
\left[\begin{array}{cc}
\varphi&0\\
0&\psi\\
\end{array}\right]\alpha-
\alpha^\ast
\left[\begin{array}{cc}
\psi&0\\
0&\varphi\\
\end{array}\right]\alpha,\end{array}\]
where $\alpha=\left[\begin{array}{cc}
\frac1{\sqrt{2}}&\frac1{\sqrt{2}}\\
\frac1{\sqrt{2}}&-\frac1{\sqrt{2}}\\
\end{array}\right]\otimes 1_n$, we get
\[\begin{array}{rcl}
Q_n(f)&=&tQ_{2n}^{'}\bigg(\frac{1}{t}\alpha^\ast\left[\begin{array}{cc}
f&0\\
0&-f\\
\end{array}\right]\alpha\bigg)\\
&=&tD_{2n}\left(\alpha^\ast
\left[\begin{array}{cc}
\varphi&0\\
0&\psi\\
\end{array}\right]\alpha,
\alpha^\ast
\left[\begin{array}{cc}
\psi&0\\
0&\varphi\\
\end{array}\right]\alpha\right)\\
&=&tD_{2n}\left(
\left[\begin{array}{cc}
\varphi&0\\
0&\psi\\
\end{array}\right],
\left[\begin{array}{cc}
\psi&0\\
0&\varphi\\
\end{array}\right]\right)\\
&=&tQ_{2n}^{'}\bigg(\frac{1}{t}\left[\begin{array}{cc}
f&0\\
0&-f\\
\end{array}\right]\bigg)\\
&=&tQ_n^{'}(\frac{1}{t} f)=Q_n^{''}(f).\end{array}\]

For $f\in U_n$ and $c\in\mathbb C$ with $c=r\exp(i\theta)$, where
$r>0$ and $\theta\in \mathbb R$, we have
\[\begin{array}{rcl}
Q_n(cf)&=&Q_{2n}^{''}\bigg(\left[\begin{array}{cc}
0&cf\\
(cf)^*&0\\
\end{array}\right]\bigg)=Q_{2n}^{''}\bigg(\gamma\left[\begin{array}{cc}
0&f\\
f^*&0\\
\end{array}\right]\gamma^*\bigg)\\
&\le&\abs{c}Q_{2n}^{''}\bigg(\left[\begin{array}{cc}
0&f\\
f^*&0\\
\end{array}\right]\bigg)=\abs{c}Q_n(f),
\end{array}\]
where $\gamma=\left[\begin{array}{cc}
\sqrt{r}\exp{(i\frac{\theta}{2})}&0\\
0&\sqrt{r}\exp{(-i\frac{\theta}{2})}
\end{array}\right]$, and hence
$Q_n(f)=Q_n(c^{-1}cf)\le\abs{c}^{-1}Q_n(cf)$. So
\[Q_n(cf)=\abs{c}Q_n(f).\]
Obviously, this holds for $c=0$. It is clear that $Q_n(f+g)\le
Q_n(f)+Q_n(g)$ for $f, g \in U_n$. Now assume that $f\in U_n, g\in
U_m$ and $\alpha\in M_{m, n}$ and $\beta\in M_{n,m}$. Then we have
\[\begin{array}{rcl}
Q_m(\alpha
f\beta)&=&\|\alpha\|\|\beta\|Q_m\left(\frac{\alpha}{\|\alpha\|}f\frac{\beta}
{\|\beta\|}\right)\\
&=&\|\alpha\|\|\beta\| Q_{2m}^{''}\Bigg(\left[\begin{array}{cc}
0&\frac{\alpha f\beta}{\|\alpha\|\|\beta\|}\\
\left(\frac{\alpha f\beta}{\|\alpha\|\|\beta\|}\right)^*&0\\
\end{array}\right]\Bigg)\\
&=&\|\alpha\|\|\beta\|Q_{2m}^{''}\Bigg( \left[\begin{array}{cc}
\frac{\alpha^*}{\|\alpha\|}&0\\
0&\frac{\beta}{\|\beta\|}
\end{array}\right]^*
\left[\begin{array}{cc}
0&f\\
f^*&0\\
\end{array}\right]
\left[\begin{array}{cc}
\frac{\alpha^*}{\|\alpha\|}&0\\
0&\frac{\beta}{\|\beta\|}\\
\end{array}\right]\Bigg)\\
&\le&\|\alpha\|\|\beta\|\Bigg\| \left[\begin{array}{cc}
\frac{\alpha^*}{\|\alpha\|}&0\\
0&\frac{\beta}{\|\beta\|}\\
\end{array}\right]\Bigg\|^2Q_{2n}^{''}\bigg(
\left[\begin{array}{cc}
0&f\\
f^*&0
\end{array}\right]\bigg)\\
&=&\|\alpha\|\|\beta\|Q_n(f)
\end{array}\] and
\[\begin{array}{rcl}
&&Q_{n+m}(f\oplus g)=Q_{2(n+m)}^{''}\bigg( \left[\begin{array}{cc}
0&f\oplus g\\
f^*\oplus g^*&0
\end{array}\right]\bigg)\\
&=&Q_{2(n+m)}^{''}\left( \left[\begin{array}{cccc}
1&0&0&0\\
0&0&1&0\\
0&1&0&0\\
0&0&0&1
\end{array}\right]\left[\begin{array}{cccc}
0&f&0&0\\
f^*&0&0&0\\
0&0&0&g\\
0&0&g^*&0
\end{array}\right]
\left[\begin{array}{cccc}
1&0&0&0\\
0&0&1&0\\
0&1&0&0\\
0&0&0&1
\end{array}\right]\right)\\
&=&Q_{2(n+m)}^{''}\left( \left[\begin{array}{cccc}
0&f&0&0\\
f^*&0&0&0\\
0&0&0&g\\
0&0&g^*&0
\end{array}\right]\right)\\
&=&\max\bigg\{Q_{2n}^{''}\bigg( \left[\begin{array}{cc}
0&f\\
f^*&0
\end{array}\right]\bigg),
Q_{2m}^{''}\bigg( \left[\begin{array}{cc}
0&g\\
g^*&0
\end{array}\right]\bigg)\bigg\}\\
&=&\max\{Q_m(f), Q_n(g)\}.
\end{array}\]
By Lemma \ref{le:41},
$U_n=M_n((\tilde{\mathcal{V}})^*)(n\in\mathbb{N})$. Therefore,
$\mathcal Q=(Q_n)$ is a matrix norm on $(\widetilde{\mathcal
V})^*$. And
\[\begin{array}{rcl}
Q_n(f^*)&=&Q_{2n}\bigg( \left[\begin{array}{cc}
0&f^*\\
f&0
\end{array}\right]\bigg)\\
&=&Q_{2n}\bigg( \left[\begin{array}{cc}
0&1\\
1&0
\end{array}\right]
\left[\begin{array}{cc}
0&f\\
f^*&0
\end{array}\right]
\left[\begin{array}{cc}
0&1\\
1&0
\end{array}\right]\bigg)\\
&=&Q_{2n}\bigg( \left[\begin{array}{cc}
0&f\\
f^*&0
\end{array}\right]\bigg)=Q_n(f),\end{array}\] for all $f\in
M_n((\tilde{\mathcal{V}})^*)$.
\end{proof}

\section{Matrix Lipschitz seminorms from matrix metrics}\label{sec:6}

The purpose of this section is to characterize the matrix metrics
on $\mathcal{CS(V)}$ which come from matrix Lipschitz seminorms.

Let $(\mathcal{V}, 1)$ be a matrix order unit space, and let
$\mathcal{D}=(D_n)$ be a matrix metric on $\mathcal{CS(V)}$. We
will refer to the topology on $\mathcal{CS(V)}$ defined by
$\mathcal{D}$, that is, the topology on each $CS_n(\mathcal{V})$
is induced by $D_n$, as the {\it $\mathcal{D}$-topology}. The
natural topology on $\mathcal{CS(V)}$ is the {\it BW-topology},
that is, topologies each $CS_n(\mathcal{V})$ by
BW-topology\cite{ar}.

Let $V$ and $W$ be two vector spaces, and let $\mathbf K=(K_n)$ be
a matrix convex set of $V$. A {\it matrix affine mapping} from
$\mathbf K$ to $W$ is a sequence $\mathbf F=(F_n)$ of mappings
$F_n:\, K_n\longmapsto M_n(W)$ such that
\[F_n\left(\sum_{i=1}^k\gamma^*_iv_i\gamma_i\right)=\sum_{i=1}^k\gamma^*_iF_{n_i}(v_i)\gamma_i,\]
for all $v_i\in K_{n_i}$ and $\gamma_i\in M_{n_i, n}$ for
$i=1,2,\cdots,k$ satisfying $\sum_{i=1}^k\gamma^*_i\gamma_i=1_n$.

A matrix convex set $\mathbf K=(K_n)$ of a locally convex vector
space $V$ is {\it compact} if each $K_n$ is compact in the product
topology in $M_n(V)$. Given a compact matrix convex set $\mathbf
K=(K_n)$ of a locally convex vector space. For $r\in\mathbb N$,
let $A(\mathbf K, M_r)$ to be the set of all matrix affine
mappings $\mathbf F=(F_n):\,\mathbf K\longmapsto M_r$, such that
$F_1$ is continuous. Using the linear structure and the adjoint
operation in $M_n(M_r)$, $A(\mathbf K, M_r)$ becomes a vector
space with a *-operation under pointwise operations. Similarly,
the order structure in $M_n(M_r)$ defines a positive cone in
$A(\mathbf K, M_r)$, where $\mathbf F\ge 0$ in $A(\mathbf K, M_r)$
if $F_n(v)\ge 0$ for all $n\in\mathbb N$ and $v\in K_n$. So
$A(\mathbf K, M_r)$ is a partially ordered space. If we denote
$\mathbf I^{(r)}=(I_n^{(r)}):\,\mathbf K\longmapsto M_r$ by
$I_n^{(r)}(v)=1_n\otimes 1_r$ for $v\in K_n$, then $\mathbf
I^{(r)}\in A(\mathbf K, M_r)$. We use the abbreviation $\mathbf
I=\mathbf I^{(1)}$. For $\mathbf F=(F_n)\in A(\mathbf K, \mathbb
C)$ and any unit vector $\xi\in{\mathbb C}^n$ and $v\in K_n$, we
have
\[<F_n(v)\xi, \xi>=\xi^*F_n(v)\xi=F_1(\xi^*v\xi).\]
So $\mathbf I$ is an order unit for $A(\mathbf K, \mathbb C)$. It
is clear that $A(\mathbf K, \mathbb C)$ is an order unit space.
Let $A(K_1, \mathbb C)$ be the space of all continuous affine
mappings from $K_1$ to $\mathbb C$. Then $A(K_1, \mathbb C)$ can
be made into an order unit space\cite{ka}.  Define $\Phi
:\,A(\mathbf K, \mathbb C)\longmapsto A(K_1, \mathbb C)$ by
$\Phi(\mathbf F)=F_1$ for $\mathbf F\in A(\mathbf K, \mathbb C)$.
Then $\Phi$ is a unital order preserving bijection of $A(\mathbf
K, \mathbb C)$ onto $A(K_1, \mathbb C)$. Identifying
$M_n(A(\mathbf K, \mathbb C))$ with $A(\mathbf K, M_n)$, we may
use the ordering on $A(\mathbf K, M_n)$ to define a positive cone
in $M_n(A(\mathbf K, \mathbb C))$. In this way $A(\mathbf K,
\mathbb C)$ becomes a matrix order unit space(see page 313 in \cite{wewi}). We will
simply denote the matrix order unit space $A(\mathbf K, \mathbb
C)$ by $A(\mathbf K)$.

A non-commutative
version of the representation theorem of Kadison was proved in \cite{wewi}(Proposition 3.5 in \cite{wewi}):

\begin{proposition}\label{pro:61}  {\rm (1)} If $\mathcal{R}$ is an operator system, then $\mathcal{CS(R)}$ is a
self-adjoint compact matrix convex set in $\mathcal R^*$, equipped
with the BW-topology, and $A(\mathcal{CS(R)})$ and
$\overline{\mathcal{R}}$, the completion of $\mathcal R$, are
isomorphic as operator systems.

{\rm (2)} If $\mathbf{K}$ is a compact matrix convex set in a
locally convex space $V$, then $A(\mathbf{K})$ is an operator
system, and $\mathbf{K}$ and $\mathcal{CS}(A(\mathbf{K}))$ are
matrix affinely homeomorphic.
\end{proposition}

So we can view the matrix order unit spaces as exactly the dense
subspaces containing the order unit $\mathbf{I}$ inside $A(\mathbf
K)$, where $\mathbf K$ is any compact matrix convex subset of a
topological vector space. The following is a direct consequence of Lemma 2.3.4 in \cite{efruan}:

\begin{lemma}\label{le:62} Suppose that $V$ is an operator space with matrix norm $\|\cdot\|=(\|\cdot\|_n)$.
Let $V^*$ be the operator space dual of $V$, with matrix norm
$\|\cdot\|^*=(\|\cdot\|_n^*)$. If $v_0\in M_n(V)$, then
\[\|v_0\|_n=\sup\{\|\ll\varphi, v_0\gg\|:\, \|\varphi\|_r^*\le 1, \varphi\in M_r(V^*),
r\in\mathbb N\}.\]
\end{lemma}

\begin{lemma}\label{le:63} Suppose that $V$ is a self-adjoint operator space with matrix norm
$\|\cdot\|=(\|\cdot\|_n)$. Let $V^*$ be the operator space dual of
$V$, with matrix norm $\|\cdot\|^*=(\|\cdot\|_n^*)$. If
$\varphi\in M_n(V^*)$, then
\[\|\varphi\|_n^*=\sup\{\|\ll\varphi, v\gg\|:\, v^*=v\in M_r(V), \|v\|_r\le 1, r\in\mathbb N\}.\]
\end{lemma}

\begin{proof}
By Proposition 2.1 in \cite{ble}, we get
\[\|\varphi\|_n^*=\sup\{\|\ll\varphi, v\gg\|:\, v\in M_r(V), \|v\|_r\le 1, r\in\mathbb N\}.\]
For $v\in M_r(V)$ with $\|v\|_r\le 1$, since
\[\left\|\left[\begin{array}{cc} 0&v\\
v^*&0\end{array}\right]\right\|_{2r}=\left\|\left[\begin{array}{cc} v&0\\
0&v^*\end{array}\right]\right\|_{2r}=\max\{\|v\|_r,
\|v^*\|_r\}=\|v\|_r\le 1,\] and
\[\begin{array}{rcl}
&&\left\|\ll\varphi, \left[\begin{array}{cc} 0&v\\
v^*&0\end{array}\right]\gg\right\|=\left\|\left[\begin{array}{cc} 0&\ll\varphi, v\gg\\
\ll\varphi, v^*\gg&0\end{array}\right]\right\|\\
&=&\left\|\left[\begin{array}{cc} \ll\varphi, v\gg&0\\
0&\ll\varphi, v^*\gg\end{array}\right]\right\|\ge\|\ll\varphi,
v\gg\|,\end{array}\] we obtain that
\[\|\varphi\|_n^*\le\sup\{\|\ll\varphi, v\gg\|:\, v^*=v\in M_r(V), \|v\|_r\le 1, r\in\mathbb N\},\]
and we are done.
\end{proof}

\begin{definition}\label{def:64} Let $(\mathcal V, 1)$ be a matrix order unit space.
By a {\it matrix Lip-norm} on $(\mathcal V, 1)$ we mean a matrix
Lipschitz seminorm $\mathcal L=(L_n)$ on $(\mathcal V, 1)$ such
that the $\mathcal{D_L}$-topology on $\mathcal{CS(V)}$ agrees with
the BW-topology.
\end{definition}

\begin{example}\label{ex:65}
The matrix Lipschitz seminorm $\mathcal{L}=(L_n)$ on
$(\mathcal{V}, 1)$ as defined in Example \ref{ex:32} is clearly
lower semicontinuous. From Proposition 3.1 in \cite{wu}, Theorem
2.3 in \cite{ri0} and its proof, we see that the image of
$L_1^{1}$ in $\tilde{\mathcal{V}}$ is totally bounded for
$\|\cdot\|_1^{\sim}$ by Theorem 1.8 in \cite{ri0}. According to
Theorem 5.3 in \cite{wu}, we see that the $\mathcal{D_L}$-topology
on $\mathcal{CS(V)}$ agrees with the BW-topology. Therefore,
$\mathcal{L}$ is a lower semicontinuous matrix Lip-norm on
$(\mathcal{V}, 1)$.
\end{example}

Let $(\mathcal{V}, 1)$ be a matrix order unit space. Let $\mathcal
A$ be a unital $C^*$-algebra such that $\mathcal V\subseteq
\mathcal A$ and $1_{\mathcal A}=1$. Fix $\rho\in CS_1(\mathcal
A)$. We define $\Phi _\rho ^{(r)}:\,\mathcal A\longmapsto M_r$ by
\[\Phi _\rho ^{(r)}(a)=\left[\begin{array}{ccc}
\rho (a)&&0\\
&\ddots&\\
0&&\rho (a)
\end{array}\right],\ \ \ a\in\mathcal A.\]
Then for $[a_{ij}]\in M_r(\mathcal A)$, we have
\[S_{\Phi _\rho ^{(r)}}([a_{ij}])=\frac1r\sum_{i=1}^r\rho (a_{ii}).\]
It is clear that $S_{\Phi _\rho ^{(r)}}$ is positive. So $\Phi
_\rho ^{(r)}$ is completely positive(Theorem 5.1 in \cite{pau})
and $\Phi _\rho ^{(r)}(1)=1_r$, that is, $\Phi _\rho ^{(r)}\in
CS_r(\mathcal A)$. Hence, $\Phi _\rho ^{(r)}\in CS_r(\mathcal V)$.

Let ${\mathcal D}=(D_n)$ be a matrix metric on ${\mathcal
C\mathcal S}(\mathcal V)$, we can define a sequence ${\mathcal
L}_{\mathcal{D}}=(L_{D_n})$ of gauges on $A({\mathcal C\mathcal
S}(\mathcal V))$ by
\[L_{D_n}(\mathbf F^{(n)})=\sup\left\{\frac{\|F^{(n)}_r(\varphi)-F^{(n)}_r(\psi)\|}{D_r(\varphi,
\psi)}:\,\varphi\neq\psi, \varphi, \psi\in CS_r(\mathcal V),
r\in\mathbb N\right\},\] where $\mathbf F^{(n)}\in
A(\mathcal{CS(V)}, M_n)$, and then define a sequence
$\mathcal{L}_{\mathcal{D}}^c=(L_{D_n}^c)$ of gauges on
$(\tilde{\mathcal V})^*$ by
\[\begin{array}{rcl}
L^c_{D_n}(g)&=&\sup\{t\|(F^{(r)}_n(\varphi_1)-F^{(r)}_n(\varphi_2))+i(F^{(r)}_n(\varphi_3)-
F^{(r)}_n(\varphi_4))\|:\\
&&L_{D_r}(\mathbf F^{(r)})\le 1, \mathbf F^{(r)}\in
A(\mathcal{CS(V)}, M_r), r\in\mathbb N\},
\end{array}\]
where $g=t[(\varphi_1-\varphi_2)+i(\varphi_3-\varphi_4)]\in
M_n((\tilde{\mathcal V})^*)$ with $\varphi_1, \varphi_2,
\varphi_3, \varphi_4\in CS_n(\mathcal V)$ and $t>0$. Since we
identify $M_n(A(\mathcal{CS}(\mathcal{V})))$ and $A( \mathcal{CS}
(\mathcal{V}), M_n)$, for $\mathbf{F}^{(n)}\in
A(\mathcal{CS}(\mathcal{V}), M_n)$ and $\gamma\in M_{n,r}$ we have
that $\gamma^*\mathbf{F}^{(n)}\gamma\in
A(\mathcal{CS}(\mathcal{V}), M_r)$ and
$(\gamma^*\mathbf{F}^{(n)}\gamma)_p(\varphi)=(\gamma\otimes
1_p)^*F^{(n)}_p(\varphi)(\gamma\otimes 1_p)$ for $\varphi\in
CS_p(\mathcal{V})$. Clearly, ${\mathcal L}_{\mathcal{D}}$ and
$\mathcal{L}_{\mathcal{D}}^c$ are matrix gauges and the null space
of each $L_{D_n}$ is $M_n(\mathbb{C}\mathbf{I})$. Denote
\[K_n=\{\mathbf F^{(n)}\in A(\mathcal{CS(V)}, M_n):\, L_{D_n}(\mathbf F^{(n)})<+\infty\},\ \ n\in
\mathbb{N},\] and $\mathbf{K}=(K_n)$.

Suppose that $\mathcal D$ is convex, midpoint balanced, and midpoint
concave. Then, by Theorem \ref{th:55}
there exists a matrix norm $\mathcal{Q}=(Q_n)$ on
$(\tilde{\mathcal{V}})^*$ such that
$Q_n(\varphi-\psi)=D_n(\varphi,\psi)$ for all $\varphi,\psi\in CS_n(\mathcal{V})$, and $Q_n(f^*)=Q_n(f)$ for
each $f\in M_n((\tilde{\mathcal{V}})^*)$ and all $n\in\mathbb{N}$.
Let $\mathcal W$ denote the operator space $(\tilde{\mathcal
V})^*$ with the matrix norm $\mathcal Q=(Q_n)$. Let $\mathcal W^*$
be the operator space dual of $\mathcal W$, with the dual operator
space norm $\mathcal Q^{d}=(Q_n^{d})$. Now we assume that the
$\mathcal{D}$-topology on $\mathcal{CS(V)}$ agrees with the
BW-topology.

\begin{proposition}\label{pro:66} \rm{(1)} $\mathbf K=(K_n)$ is an absolutely matrix convex set in
$A(\mathcal{CS(V)})$.

\rm{(2)} for $g\in M_n(\mathcal W^*)$, define $\sigma^{(n)}(g):\,
\mathcal{CS(V)}\longmapsto M_n$ by
\[\sigma^{(n)}_r(g)(\varphi)=\ll g, \varphi-\Phi _\rho ^{(r)}\gg,\]
for $\varphi\in CS_r(\mathcal V)$ and $r\in\mathbb{N}$, then
$\sigma^{(n)}(g)\in K_n$ and $L_{D_n}(\sigma^{(n)}(g))\le Q_n^d(g)$.

\rm{(3)} $\overline{K}_n=A(\mathcal{CS(V)}, M_n)$, where
$\overline{\mathbf K}=(\overline{K}_n)$ is the closure of $\mathbf
K$.
\end{proposition}

\begin{proof} (1) This follows directly from the fact that $\mathcal{L}_{\mathcal{D}}$ is a matrix gauge on 
$A(\mathcal{CS(V)})$.

(2) It is easy to see that $M_n(\mathbb C\mathbf I)\subseteq K_n$.
Let $\varphi_i\in CS_{m_i}(\mathcal V)$ and $\gamma_i\in M_{m_i,
r},\, i=1, 2, \cdots, k$, satisfying $\sum_{i=1}^k
\gamma^*_i\gamma_i=1_r$. Then we have
\[\begin{array}{rcl}
\sigma_r^{(n)}(g)(\sum_{i=1}^k\gamma_i^*\varphi_i\gamma_i)
&=&\ll g, \sum_{i=1}^k\gamma^*_i\varphi_i\gamma_i-\Phi _\rho ^{(r)}\gg\\
&=&\ll g, \sum_{i=1}^k\gamma^*_i(\varphi_i-\Phi _\rho ^{(m_i)})\gamma_i\gg\\
&=&\sum_{i=1}^k (\gamma_i\otimes 1_n)^*\ll g, \varphi_i-\Phi
_\rho^{(m_i)}\gg(\gamma_i
\otimes 1_n)\\
&=&\sum_{i=1}^k (\gamma_i\otimes 1_n)^*
\sigma_{m_i}^{(n)}(g)(\varphi_i)(\gamma_i\otimes 1_n),
\end{array}\]
that is, $\sigma^{(n)}(g)$ is a matrix affine mapping from
$\mathcal{CS(V)}$ to $M_n$. For $\varphi, \psi\in CS_r(\mathcal
V)$, we have
\[\begin{array}{rcl}
&&\|\sigma^{(n)}_r(g)(\varphi)-\sigma^{(n)}_r(g)(\psi)\|=\|\ll g, \varphi-\psi\gg\|\\
&\le& Q^d_n(g)Q_r(\varphi-\psi)=Q_n^d(g)D_r(\varphi,
\psi).\end{array}\] by Lemma \ref{le:62}. So $\sigma^{(n)}_1(g)$
is continuous on $CS_1(\mathcal V)$ since $\mathcal D$ gives the
BW-topology. Thus $\sigma^{(n)}(g)\in A(\mathcal{CS(V)}, M_n)$.
Moreover, we have that $\sigma^{(n)}(g)\in K_n$ and
$L_{D_n}(\sigma^{(n)}$ $(g))\le Q^d_n(g)$.

(3) Let $\overline{\mathbf K}=(\overline{K}_n)$ be the closure of
$\mathbf K$. Assume that $f\in [A(\mathcal{CS(V)}, M_n)/M_n(\mathbb{C}\mathbf{I})]$
$\setminus[\overline{K}_n/M_n (\mathbb{C}\mathbf{I})]$. By the matricial separation theorem(Theorem 1.6 in
\cite{wewi}), applied to $A(\mathcal{CS(V)})/(\mathbb C\mathbf I)$
and the closed matrix convex set $\tilde{\overline{\mathbf
K}}=(\overline{K}_n/M_n(\mathbb C\mathbf I))$ in
$A(\mathcal{CS(V)})/(\mathbb C\mathbf I)$, there is a continuous
linear mapping $\Phi :\, A(\mathcal{CS(V)})/(\mathbb C\mathbf
I)\longmapsto M_n$ such that
\[\mathrm{Re}\Phi_r(g)\le 1_n\otimes 1_r\]
for all $r\in\mathbb N$, $g\in\overline{K}_r/M_r(\mathbb C\mathbf
I)$ and
\[\mathrm{Re}\Phi_n(f)\not\leq 1_n\otimes 1_n.\]
Since $\sigma^{(r)}(M_r(\mathcal W^*))\subseteq K_r$ for $r\in
\mathbb N$, we have
\[\mathrm{Re}\Phi_r(\sigma^{(r)}(g)+\mathbb{C}\mathbf{I}^{(r)})\le 1_n\otimes 1_r\]
for all $r\in\mathbb N$, $g\in M_r(\mathcal W^*)$. Identifying
$\Phi$ with $F\in M_n((A(\mathcal{CS(V)})/(\mathbb C\mathbf
I))^*)\cong M_n((\tilde{\mathcal V})^*)= M_n(\mathcal W)$ by
Proposition \ref{pro:61}, this means that
\[\mathrm{Re}\ll F, g\gg\le 1_n\otimes 1_r.\]
So by the operator space duality of $\mathcal W$ and $\mathcal
W^*$, we obtain $F=0_n$, that is, $\Phi(g)=0$ for $g\in
A(\mathcal{CS(V)})/(\mathbb C\mathbf I)$. Therefore,
\[\mathrm{Re}\Phi_n(f)=0_{n\times n}\le 1_n\otimes 1_n,\]
a contradiction. Since $M_n(\mathbb C\mathbf I)\subseteq
A(\mathcal{CS(V)}, M_n) \cap\overline{K}_n$,
$\overline{K}_n=A(\mathcal{CS(V)}, M_n)$.
\end{proof}

\begin{proposition}\label{pro:67} If
$L^c_{D_n}(\varphi-\psi)=D_n(\varphi, \psi)$ for $\varphi, \psi\in
CS_n(\mathcal{V})$, then $\mathcal{L}_{\mathcal{D}}$ is a lower
semicontinuous matrix Lip-norm on $(K_1, I_1^{(1)})$.
\end{proposition}

\begin{proof}
It is clear that $(K_1, I_1^{(1)})$ is a matrix order unit space
as a subspace of $A(\mathcal{CS(V)})$ and
$\mathcal{L}_{\mathcal{D}}=(L_{D_{n}})$ is a matrix Lipschitz
seminorm on it. Similar to the proof of
Proposition \ref{pro:36}, we have that
$\mathcal{L}_{\mathcal{D}}=(L_{D_n})$ is lower semicontinuous.
Since $K_1$ is dense in $A(\mathcal{CS}(\mathcal{V}))$,
$\mathcal{CS}({A}(\mathcal{CS(V)}))=\mathcal{CS}(K_1)$. 
By Proposition \ref{pro:61}(2), we obtain that
$\mathcal{CS}(K_1)=\mathcal{CS(V)}$. Since
$L^c_{D_n}(\varphi-\psi)=D_n(\varphi, \psi)$ for $\varphi, \psi\in
CS_n(\mathcal{V})$ and clearly
$L^c_{D_n}(\varphi-\psi)=D_{L_{D_n}}(\varphi,\psi)$ for $\varphi,
\psi\in CS_n(\mathcal{V})$, the
$\mathcal{D}_{\mathcal{L}_{\mathcal{D}}}$-topology on
$\mathcal{CS}(K_1)$ agrees with the BW-topology. Therefore,
$\mathcal{L}_{\mathcal{D}}$ is a matrix Lip-norm on $(K_1,
I_1^{(1)})$.
\end{proof}

For $\mathbf F^{(n)}=(F_r^{(n)})\in K_n$, define $\rho ^{(n)}_{r,
1}(\mathbf F^{(n)}): \, SB^2_r(\tilde{\mathcal{V}})\longmapsto
M_r(M_n)$ by
\[\rho ^{(n)}_{r, 1}(\mathbf F^{(n)})(\varphi)=F_r^{(n)}(\varphi_1)-F_r^{(n)}(\varphi_2),\]
where $\varphi=\varphi_1-\varphi_2\in SB^2_r(\tilde{\mathcal{V}})$
and $\varphi_1, \varphi_2\in CS_r(\mathcal V)$. If also
$\varphi=\varphi_3-\varphi_4$ for some $\varphi_3, \varphi_4\in
CS_r(\mathcal V)$, then
$\frac12(\varphi_1+\varphi_4)=\frac12(\varphi_2+\varphi_3)$. These
are elements of $CS_r(\mathcal V)$ and so
$F_r^{(n)}\left(\frac12(\varphi_1+\varphi_4)\right)=F_r^{(n)}\left(\frac12(\varphi_2+\varphi_3)\right)$.
From the fact that $\mathbf F^{(n)}$ is matrix affine it now
follows that
$F_r^{(n)}(\varphi_1)-F_r^{(n)}(\varphi_2)=F_r^{(n)}(\varphi_3)-F_r^{(n)}(\varphi_4)$.
Thus $\rho^{(n)}_*(\mathbf F^{(n)})=(\rho^{(n)}_{r, 1}(\mathbf
F^{(n)}))$ is well-defined.

For $\varphi=\varphi_1-\varphi_2\in SB^2_r(\tilde{\mathcal{V}})$
with $\varphi_1, \varphi_2\in CS_r (\mathcal V)$, we have
\[\|\rho^{(n)}_{r, 1}(\mathbf
F^{(n)})(\varphi)\|=\|F_r^{(n)}(\varphi_1)-F_r^{(n)}
(\varphi_2)\|\le L_{D_n}(\mathbf{F}^{(n)})Q_r(\varphi).\] If $t\in
[-1, 1]$ and $f=\varphi-\psi\in SB^2_r(\tilde{\mathcal{V}})$ with
$\varphi, \psi\in CS_r(\mathcal V)$, we have that $
tf=\mathrm{sgn}(t)(\varphi-(\abs{t}\psi+(1-\abs{t})\varphi))\in
SB^2_r(\tilde{\mathcal{V}})$ and
\[\begin{array}{rcl}
\rho^{(n)}_{r, 1}(\mathbf F^{(n)})(tf)&=&\rho^{(n)}_{r, 1}(\mathbf F^{(n)})
(\mathrm{sgn}(t)(\varphi-(\abs{t}\psi+(1-\abs{t})\varphi)))\\
&=&\mathrm{sgn}(t)[F_r^{(n)}(\varphi)-F_r^{(n)}(\abs{t}\psi+(1-\abs{t})\varphi)]\\
&=&\mathrm{sgn}(t)[F_r^{(n)}(\varphi)-\abs{t}F_r^{(n)}(\psi)-(1-\abs{t})F_r^{(n)}(\varphi)]\\
&=&\mathrm{sgn}(t)\abs{t}(F_r^{(n)}(\varphi)-F_r^{(n)}(\psi))=t\rho^{(n)}_{r,
1}(\mathbf F^{(n)})(f).\end{array}\] Let $f=\varphi_1-\psi_1$ and
$g=\varphi_2-\psi_2\in SB_r^2(\tilde{\mathcal{V}})$ with
$\varphi_1, \varphi_2, \psi_1, \psi_2\in CS_r(\mathcal V)$ and
$f+g\in SB^2_r(\tilde{\mathcal{V}})$. Then we have
\[\begin{array}{rcl}
\rho^{(n)}_{r, 1}(\mathbf F^{(n)})(f+g)&=&2\rho^{(n)}_{r, 1}(\mathbf F^{(n)})(\frac{f+g}2)\\
&=&2(F_r^{(n)}(\frac{\varphi_1+\varphi_2}2)-F_r^{(n)}(\frac{\psi_1+\psi_2}2))\\
&=&F_r^{(n)}(\varphi_1)+F_r^{(n)}(\varphi_2)-F_r^{(n)}(\psi_1)-F_r^{(n)}(\psi_2)\\
&=&\rho^{(n)}_{r, 1}(\mathbf F^{(n)})(f)+\rho^{(n)}_{r, 1}(\mathbf
F^{(n)})(g).
\end{array}\]

If $f\in\mathbb{R}SB^2_r(\tilde{\mathcal{V}})$, then there is a
$t>0$ such that $tf\in SB^2_r(\tilde{\mathcal{V}})$. Set
\[\rho^{(n)}_{r, 2}(\mathbf F^{(n)})(f)=t^{-1}\rho^{(n)}_{r, 1}(\mathbf F^{(n)})(tf).\]
Then it is clear that $\rho_{**}^{(n)}(\mathbf
F^{(n)})=(\rho^{(n)}_{r, 2}(\mathbf F^{(n)}))$ is well-defined and
\[\rho^{(n)}_{r, 2}(\mathbf F^{(n)})(sf)=s\rho^{(n)}_{r, 2}(\mathbf F^{(n)})(f),\]
for $s\in\mathbb{R}$ and
$f\in\mathbb{R}SB^2_r(\tilde{\mathcal{V}})$. If
$f=t_1(\varphi_1-\psi_1)\in\mathbb{R}SB^2_r(\tilde{\mathcal{V}})$
and
$g=t_2(\varphi_2-\psi_2)\in\mathbb{R}SB^2_r(\tilde{\mathcal{V}})$
with $\varphi_1, \varphi_2, \psi_1, \psi_2\in CS_r(\mathcal{V})$
and $t_1, t_2>0$, we may assume that $t_1, t_2>1$, and then we
have
\[\begin{array}{rcl}
\rho^{(n)}_{r, 2}(\mathbf F^{(n)})(f+g)&=&\rho^{(n)}_{r, 2}(\mathbf F^{(n)})
((t_1\varphi_1+t_2\varphi_2)-(t_1\psi_1+t_2\psi_2))\\
&=&\rho^{(n)}_{r, 2}(\mathbf
F^{(n)})\left((t_1+t_2)\left(\frac{t_1\varphi_1+t_2\varphi_2}{t_1+t_2}-
\frac{t_1\psi_1+t_2\psi_2}{t_1+t_2}\right)\right)\\
&=&(t_1+t_2)\rho^{(n)}_{r, 2}(\mathbf F^{(n)})\left(\frac1{t_1+t_2}f+\frac1{t_1+t_2}g\right)\\
&=&(t_1+t_2)\rho^{(n)}_{r, 2}(\mathbf
F^{(n)})\left(\frac1{t_1+t_2}f\right)+(t_1+t_2)\rho^{(n)}_{r, 2}
(\mathbf F^{(n)})\left(\frac1{t_1+t_2}g\right)\\
&=&\rho^{(n)}_{r, 2}(\mathbf F^{(n)})(f)+\rho^{(n)}_{r, 2}(\mathbf F^{(n)})(g).
\end{array}\]
For $f=tf_1\in\mathbb{R}SB^2_r(\tilde{\mathcal{V}})$ with $t\ge 1$
and $f_1\in SB^2_r(\tilde{\mathcal{V}})$, we have
\[\|\rho^{(n)}_{r,2}(\mathbf F^{(n)})(f)\|=\|t\rho^{(n)}_{r,1}(\mathbf F^{(n)})(f_1)\|
\le tL_{D_n}(\mathbf F^{(n)})Q_r(f_1)=L_{D_n}(\mathbf
F^{(n)})Q_r(f).\]

Now for
$f=f_1+if_2\in\mathbb{R}SB^2_r(\tilde{\mathcal{V}})+i\mathbb{R}SB^2_r(\tilde{\mathcal{V}})$
with $f_1, f_2\in\mathbb{R}SB^2_r(\tilde{\mathcal{V}})$, we define
\[\rho^{(n)}_r(\mathbf F^{(n)})(f)=\rho^{(n)}_{r,2}(\mathbf F^{(n)})(f_1)+i\rho^{(n)}_{r,2}
(\mathbf F^{(n)})(f_2).\] Obviously, $\rho^{(n)}(\mathbf
F^{(n)})=(\rho^{(n)}_r(\mathbf F^{(n)})$ is well-defined.

\begin{proposition}\label{pro:68} $\rho_1^{(n)}(\mathbf{F}^{(n)})\in
\mathcal{CB}(\mathcal{W}, M_n)=M_n(\mathcal{W}^{*})$.
\end{proposition}

\begin{proof} If $\alpha =\alpha _1+i\alpha _2$ with $\alpha _1, \alpha
_2\in\mathbb{R}$ and $f=f_1+if_2$, $g=g_1+ig_2
\in\mathbb{R}SB^2_r(\tilde{\mathcal{V}})+i\mathbb{R}SB^2_r(\tilde{\mathcal{V}})$
with $f_1, f_2, g_1, g_2\in\mathbb{R}SB^2_r(\tilde{\mathcal{V}})$,
we have
\[\begin{array}{rcl}
&&\rho _r^{(n)}(\mathbf F^{(n)})(\alpha f+g)\\
&=&\rho _{r,2}^{(n)}(\mathbf F^{(n)})(\alpha _1f_1-\alpha _2f_2+g_1)+i\rho _{r,2}^{(n)}
(\mathbf F^{(n)})(\alpha _1f_2+\alpha _2f_1+g_2)\\
&=&\alpha \rho _r^{(n)}(\mathbf F^{(n)})(f)+\rho _r^{(n)}(\mathbf
F^{(n)})(g).\end{array}\] For $\varphi _1, \varphi _2\in
CS_r(\mathcal V)$ and  $\alpha \in M_{r,m}$ with $\alpha ^*\alpha
\le 1_m$, denoting $\gamma=\left[\begin{array}{c}
\alpha \\
\sqrt{1_m-\alpha^*\alpha}\end{array}\right]$, we have
\[\alpha
^*(\varphi_1-\varphi_2)\alpha=\gamma^*\left[\begin{array}{cc}
\varphi_1-\varphi_2&0\\
0&0_m\end{array}\right]\gamma=\gamma^*\left[\begin{array}{cc}
\varphi_1&0\\
0&\Phi_\rho^{(m)}\end{array}\right]\gamma-\gamma^*\left[\begin{array}{cc}
\varphi_2&0\\
0&\Phi_\rho^{(m)}\end{array}\right]\gamma.\] Then
\[\begin{array}{rcl}
&&\rho_m^{(n)}(\mathbf F^{(n)})(\alpha^*(\varphi _1-\varphi_2)\alpha)\\
&=&F^{(n)}_m\left(\gamma^*\left[\begin{array}{cc}
\varphi_1&0\\
0&\Phi_\rho^{(m)}\end{array}\right]\gamma\right)-F_m^{(n)}\left(\gamma^*\left[\begin{array}{cc}
\varphi_2&0\\
0&\Phi_\rho^{(m)}\end{array}\right]\gamma\right)\\
&=&(\gamma\otimes 1_n)^*F^{(n)}_{m+r}\left(\left[\begin{array}{cc}
\varphi_1&0\\
0&\Phi_\rho^{(m)}\end{array}\right]\right)(\gamma\otimes 1_n)\\
&&-(\gamma\otimes 1_n)^*F_{m+r}^{(n)}\left(\left[\begin{array}{cc}
\varphi_2&0\\
0&\Phi_\rho^{(m)}\end{array}\right]\right)(\gamma\otimes 1_n)\\
&=&(\gamma\otimes 1_n)^*F^{(n)}_{m+r}\left(\xi_1
\varphi_1\xi_1^*+\xi_2\Phi_\rho ^{(m)}\xi_2^*\right)(\gamma\otimes 1_n)\\
&&-(\gamma\otimes
1_n)^*F^{(n)}_{m+r}\left(\xi_1\varphi_2\xi_1^*+\xi_2\Phi_\rho
^{(m)}\xi_2^*\right)(\gamma\otimes 1_n)\\
&=&(\gamma\otimes 1_n)^*\Big((\xi_1\otimes 1_n)
F^{(n)}_{r}(\varphi_1)(\xi_1\otimes 1_n)^*\\
&&+(\xi_2\otimes 1_n)F_m^{(n)}(\Phi _\rho ^{(m)})(\xi_2\otimes
1_n)^*\Big)(\gamma\otimes 1_n)\\
&&-(\gamma\otimes 1_n)^*\Big((\xi_1\otimes
1_n)F_r^{(n)}(\varphi_2)(\xi_1\otimes 1_n)^*\\
&&+(\xi_2\otimes 1_n)F_m^{(n)}(\Phi _\rho
^{(m)})(\xi_2\otimes 1_n)^*\Big)(\gamma\otimes 1_n)\\
&=&(\alpha\otimes 1_n)^*F_r^{(n)}(\varphi_1)(\alpha\otimes 1_n) -(\alpha\otimes 1_n)^*F_r^{(n)}(\varphi_2)(\alpha\otimes 1_n) \\
&=&(\alpha\otimes 1_n)^*\rho_r^{(n)}(\mathbf
F^{(n)})(\varphi_1-\varphi_2)(\alpha\otimes 1_n),\end{array}\]
where $\xi_1=\left[\begin{array}{c}
1_r\\
0_{m,r}\end{array}\right]$ and $\xi_2=\left[\begin{array}{c}
0_{r,m}\\
1_m\end{array}\right]$. So for $\varphi\in
SB^2_r(\tilde{\mathcal{V}})$ and $\alpha\in M_{r,m}$, we have
\[\begin{array}{rcl}
\rho_m^{(n)}(\mathbf
F^{(n)})(\alpha^*\varphi\alpha)&=&(1+\|\alpha\|)^2\rho_m^{(n)}(\mathbf
F^{(n)})\left(\frac{\alpha^*}{\|\alpha\|+1}\varphi\frac{\alpha}{\|\alpha\|+1}\right)\\
&=&(\alpha\otimes 1_n)^* \rho_r^{(n)}(\mathbf
F^{(n)})(\varphi)(\alpha\otimes 1_n).\end{array}\] Let $\varphi\in
SB^2_r(\tilde{\mathcal{V}})$ and $\alpha\in M_{m,r}$, $\beta\in
M_{r,m}$. Since
\[\begin{array}{rcl}
\alpha\varphi\beta&=&\frac12((\alpha^*+\beta)^*\varphi (\alpha ^*+\beta )-i(\alpha^*
+i\beta)^*\varphi(\alpha ^*+i\beta )\\
&&-\alpha\varphi\alpha ^*-\beta ^*\varphi \beta +i\alpha \varphi
\alpha ^*+i\beta ^* \varphi \beta ),\end{array}\] we have
\[\begin{array}{rcl}
&&\rho _m^{(n)}(\mathbf F^{(n)})(\alpha\varphi\beta)\\
&=&\frac12(((\alpha^*+\beta)\otimes 1_n)^*\rho_r^{(n)}(\mathbf
F^{(n)})(\varphi)((\alpha ^*+\beta )\otimes 1_n)\\
&&-i((\alpha^*+i\beta)\otimes 1_n)^*\rho_r^{(n)}(\mathbf F^{(n)})(\varphi)((\alpha ^*+i\beta)\otimes 1_n)\\
&&-(\alpha\otimes 1_n)\rho_r^{(n)}(\mathbf
F^{(n)})(\varphi)(\alpha\otimes 1_n)^*-(\beta\otimes
1_n)^*\rho_r^{(n)}
(\mathbf F^{(n)})(\varphi)(\beta\otimes 1_n)\\
&&+i(\alpha\otimes 1_n)\rho_r^{(n)}(\mathbf
F^{(n)})(\varphi)(\alpha\otimes 1_n)^*+i(\beta\otimes
1_n)^*\rho_r^{(n)}
(\mathbf F^{(n)})(\varphi)(\beta\otimes 1_n))\\
&=&(\alpha\otimes 1_n)\rho_r^{(n)}(\mathbf
F^{(n)})(\varphi)(\beta\otimes 1_n).\end{array}\] So for
$f=f_1+if_2\in\mathbb{R}SB^2_r(\tilde{\mathcal{V}})+i\mathbb{R}SB^2_r(\tilde{\mathcal{V}})$
and $\alpha\in M_{m,r},\,\beta \in M_{r,m}$, we have that $\alpha
f\beta =\alpha f_1\beta +i\alpha f_2\beta$, and hence
\[\begin{array}{rcl}
&&\rho_m^{(n)}(\mathbf F^{(n)})(\alpha f\beta )\\
&=&(\alpha\otimes
1_n)\rho_r^{(n)}(\mathbf F^{(n)})
(f_1)(\beta\otimes 1_n)+i(\alpha\otimes 1_n)\rho_r^{(n)}(\mathbf F^{(n)})(f_2)(\beta\otimes 1_n)\\
&=&(\alpha\otimes 1_n)\rho_r^{(n)}(\mathbf
F^{(n)})(f)(\beta\otimes 1_n).\end{array}\] Therefore, for
$f=[f_{ij}]\in M_r((\tilde{\mathcal{V}})^*)=
\mathbb{R}SB^2_r(\tilde{\mathcal{V}})+i\mathbb{R}SB^2_r(\tilde{\mathcal{V}})$,
we obtain
\[\begin{array}{rcl}
&&\rho_r^{(n)}(\mathbf F^{(n)})(f)=\rho_r^{(n)}(\mathbf
F^{(n)})\left(
\sum_{i,j=1}^re_i^*f_{ij}e_j\right)\\
&=&\sum_{i,j=1}^r\rho_r^{(n)}(\mathbf F^{(n)})(e_i^*f_{ij}e_j)
=\sum_{i,j=1}^r(e_i\otimes 1_n)^*\rho_1^{(n)}(\mathbf F^{(n)})(f_{ij})(e_j\otimes 1_n)\\
&=&[\rho_1^{(n)}(\mathbf F^{(n)})(f_{ij})]
=\ll\rho_1^{(n)}(\mathbf F^{(n)}), [f_{ij}]\gg
=\ll\rho_1^{(n)}(\mathbf F^{(n)}), f\gg ,\end{array}\] where $e_i$
is the $1\times r$ matrix $[a_{1,k}]$ with $a_{1,i}=1$ and
$a_{1,k}=0$ for $k\neq i$. Now we have
\[\begin{array}{rcl}
&&\|\ll\rho_1^{(n)}(\mathbf{F}^{(n)}),
f\gg\|=\|\rho_{r}^{(n)}(\mathbf{F}^{(n)})(f)\|\\
&\le&L_{D_n}(\mathbf{F}^{(n)})(Q_r(\mathrm{Re}f)+Q_r(\mathrm{Im}f))
\le 2L_{D_n}(\mathbf{F}^{(n)})Q_r(f)
\end{array}\] for $f\in M_r(\mathcal{W})$. It follows that $\rho_1^{(n)}(\mathbf{F}^{(n)})\in
\mathcal{CB}(\mathcal{W}, M_n)=M_n(\mathcal{W}^{*})$.\end{proof}

By Lemma \ref{le:63}, we have
\[\begin{array}{rcl}
&&Q_n^d(\rho^{(n)}_1(\mathbf F^{(n)}))\\
&=&\sup\{\|\ll \rho^{(n)}_1(\mathbf F^{(n)}), f\gg\|:\, Q_r(f)\le
1, f^*=f\in M_r(\mathcal W),
r\in\mathbb N\}\\
&=&\sup\{\|\rho^{(n)}_r(\mathbf F^{(n)})(f)\|:\, Q_r(f)\le 1,
f^*=f\in M_r(\mathcal W),
r\in\mathbb N\}\\
&=&\sup\{\|\rho^{(n)}_{r,2}(\mathbf F^{(n)})(f)\|:\, Q_r(f)\le 1,
f^*=f\in
M_r(\mathcal W), r\in\mathbb N\}\\
&\le&L_{D_n}(\mathbf F^{(n)}).\end{array}\] Let
$\tilde{K}_n=K_n/M_n(\mathbb C\mathbf I)$. Since $L_{D_n}(\mathbf
F^{(n)})=0$ for $\mathbf F^{(n)}\in M_n(\mathbb C\mathbf I)$,
$\tilde{L}_{D_n}(\tilde{\mathbf F}^{(n)})= L_{D_n}(\mathbf
F^{(n)})$ is well-defined for $\tilde{\mathbf F}^{(n)}\in
\tilde{K}_n$. $\rho_1^{(n)}(M_n(\mathbb C\mathbf I))=0_n$ implies
that $\rho_1^{(n)}$ determines a linear mapping
$\tilde{\rho}_1^{(n)}$ from $(\tilde{K}_n, \tilde{L}_{D_n})$ to
$(M_n(\mathcal W^*), Q^d_n)$ and
$Q_n^d(\tilde{\rho}_1^{(n)}(\tilde{\mathbf F}^{(n)}))
\le\tilde{L}_{D_n}(\tilde{\mathbf F}^{(n)})$.

Let $\sigma^{(n)}$ be the mapping, defined in Proposition
\ref{pro:66}, from $M_n(\mathcal{W}^*)$ to
$A(\mathcal{CS(V)}$, $M_n)$. Let $\tilde{\sigma}^{(n)}$ denote
$\sigma^{(n)}$ composed with the mapping from $A(\mathcal{CS(V)},
M_n)$ to $A(\mathcal{CS(V)}$, $M_n) /M_n(\mathbb C\mathbf I)$.
Then $\tilde{L}_{D_n}(\tilde{\sigma}^{(n)}(g))\le Q_n^d(g)$ for
$g\in M_n(\mathcal{W}^*)$. For $\mathbf F^{(n)}=(F_r^{(n)})\in
K_n$, $\varphi\in CS_r(\mathcal V)$, we have
\[\begin{array}{rcl}
&&\sigma^{(n)}_r(\rho^{(n)}_1(\mathbf
F^{(n)}))(\varphi)=\ll\rho_1^{(n)}(\mathbf F^{(n)}),
\varphi-\Phi_\rho^{(r)}\gg\\
&=&\rho_r^{(n)}(\mathbf F^{(n)})(\varphi-\Phi_\rho^{(r)})
=F_r^{(n)}(\varphi)-F_r^{(n)}(\Phi_\rho^{(r)}).\end{array}\]
Consequently,
$\tilde{\sigma}^{(n)}(\tilde{\rho}^{(n)}_1(\tilde{\mathbf
F}^{(n)}))= \tilde{\mathbf F}^{(n)}$. Similarly, for
$g=t[(\varphi_1-\varphi_2)+i(\varphi_3-\varphi_4)]\in M_n(\mathcal
W^*)$ with $t\in\mathbb R$ and $\varphi_1, \varphi_2, \varphi_3,
\varphi_4\in CS_n(\mathcal V)$ and $f\in M_r(\mathcal W)$ we have
\[\begin{array}{rcl}
&&\ll\tilde{\rho}_1^{(n)}(\tilde{\sigma}^{(n)}(g)), f\gg=\tilde{\rho}_r^{(n)}(\tilde{\sigma}^{(n)}(g))(f)\\
&=&t\left[\left(\sigma^{(n)}_r(g)(\varphi_1)-\sigma^{(n)}_r(g)(\varphi_2)\right)+i\left(\sigma^{(n)}_r(g)
(\varphi_3)-\sigma^{(n)}_r(g)(\varphi_4)\right)\right]\\
&=&t\ll g, \varphi_1-\Phi_\rho^{(r)}\gg-t\ll g, \varphi_2-\Phi_\rho^{(r)}\gg\\
&&+it\ll g, \varphi_3-\Phi_\rho^{(r)}\gg-it\ll g, \varphi_4-\Phi_\rho^{(r)}\gg\\
&=&\ll g, f\gg,\end{array}\] so that
$\tilde{\rho}_1^{(n)}(\tilde{\sigma}^{(n)}(g))=g$. Therefore,

\begin{proposition}\label{pro:69}
$\tilde{\sigma}^{(n)}$ is an isometric isomorphism of
$(M_n(\mathcal{W}^\ast), Q^d_n)$ onto $(\tilde{K}_n$, $\tilde{L}_{D_n})$ with inverse
$\tilde{\rho}^{(n)}$.\end{proposition}

\begin{theorem}\label{th:610} $\mathcal{L}_{\mathcal{D}}=(L_{D_n})$ is a lower semicontinuous
matrix Lip-norm on $(K_1, I_1^{(1)})$, and
\[D_n(\varphi, \psi)=L_{D_n}^c(\varphi-\psi)\]
for $\varphi,\psi\in CS_n(\mathcal{V})$ and $n\in\mathbb{N}$.
\end{theorem}

\begin{proof} For
$\mathbf{F}^{(n)}\in K_n$ and $\mathbf{G}^{(n)}\in
M_n(\mathbb{C}\mathbf{I})$, we have
$L_{D_n}(\mathbf{F}^{(n)}+\mathbf{G}^{(n)})=L_{D_n}(\mathbf{F}^{(n)})$,
and hence
$\tilde{L}_{D_n}(\tilde{\mathbf{F}}^{(n)})=L_{D_n}(\mathbf{F}^{(n)})$.
Now, for $g=t[(\varphi_1-\varphi_2)+i(\varphi_3-\varphi_4)]\in
M_n((\tilde{\mathcal V})^*)$ with $\varphi_1, \varphi_2,
\varphi_3, \varphi_4\in CS_n(\mathcal V)$ and $t>0$, we have
\[\begin{array}{rcl}
L_{D_n}^c(g)&=&\sup\{\|
t[(F^{(r)}_n(\varphi_1)-F^{(r)}_n(\varphi_2))+i(F^{(r)}_n(\varphi_3)-
F^{(r)}_n(\varphi_4))]\|:\\
&&L_{D_r}(\mathbf F^{(r)})\le 1, \mathbf F^{(r)}\in M_r(A(\mathcal{CS(V)})), r\in\mathbb N\}\\
&=&\sup\{\|\rho^{(r)}_n(\mathbf F^{(r)})(g)\|:\, L_{D_r}(\mathbf
F^{(r)})\le 1, \mathbf F^{(r)}
\in A(\mathcal{CS(V)}, M_r), r\in\mathbb N\}\\
&=&\sup\{\|\rho^{(r)}_n(\mathbf F^{(r)})(g)\|:\, L_{D_r}(\mathbf
F^{(r)})\le 1, \mathbf F^{(r)}
\in K_r, r\in\mathbb N\}\\
&=&\sup\{\|\tilde{\rho}^{(r)}_n(\tilde{\mathbf F}^{(r)})(g)\|:\,
\tilde{L}_{D_r}(\tilde{\mathbf
F}^{(r)})\le 1, \tilde{\mathbf F}^{(r)}\in\tilde{K}_r, r\in\mathbb N\}\\
&=&\sup\{\|\tilde{\rho}^{(r)}_n(\tilde{\sigma}^{(r)}(h))(g)\|:\,
Q^d_r(h)\le 1, h\in M_r(\mathcal
W^*), r\in\mathbb N\}\\
&=&\sup\{\|\ll h, g\gg\|:\, Q^d_r(h)\le 1, h\in M_r(\mathcal W^*), r\in\mathbb N\}\\
&=&Q_n(g)
\end{array}\]
by Proposition \ref{pro:66}. So for $\varphi, \psi\in CS_n(\mathcal
V)$, we have
\[L_{D_n}^c(\varphi-\psi)=Q_n(\varphi-\psi)=D_n(\varphi, \psi),\]
and thus by proposition \ref{pro:67},
$\mathcal{L}_{\mathcal{D}}=(L_{D_n})$ is a lower semicontinuous
matrix Lip-norm on $(K_1, I_1^{(1)})$.\end{proof}

\begin{definition}\label{de:611}By a {\it matrix Lipschitz gauge}
on a complete matrix order unit space $(\mathcal V, 1)$ we mean a
matrix gauge $\mathcal{G}=(G_n)$ on $\mathcal{V}$ such that
\begin{enumerate}
\item the null space of each $G_n$ is $M_n(\mathbb{C}1)$;
\item $G_n(v^*)=G_n(v)$ for any $v\in M_n(\mathcal{V})$;
\item $\{v\in\mathcal{V}: G_1(v)<+\infty\}$ is dense in
$\mathcal{V}$.
\end{enumerate}
A {\it matrix Lip-gauge} can be defined similarly.
\end{definition}

If $\mathcal{G}=(G_n)$ on $\mathcal{V}$ is a matrix Lipschitz
gauge on a complete matrix order unit space $(\mathcal V, 1)$, we
can similarly define a matrix metric
$\mathcal{D}_\mathcal{G}=(D_{G_n})$ on $\mathcal{CS(V)}$ by:
\[D_{G_n}(\varphi, \psi)=\sup\{\|\ll \varphi, a\gg -\ll\psi, a\gg\|:\, a\in M_r(\mathcal V),
G_r(a)\le 1, r\in{\mathbb N}\},\] where $\varphi , \psi\in
CS_n(\mathcal V)$(They also may take value $+\infty$). Given, on
the other hand, a matrix metric $\mathcal{D}=(D_n)$ on
$\mathcal{CS(V)}$, we can also define a matrix gauge
$\mathcal{G}_\mathcal{D}=(G_{D_n})$ on $(\mathcal V, 1)$ by
\[G_{D_n}(a)=\sup\left\{\frac{\|\ll \varphi, a\gg -\ll\psi, a\gg\|}{D_r(\varphi,
\psi)}:\,\varphi, \psi\in CS_r(\mathcal V), \varphi\neq \psi,
r\in{\mathbb N}\right\},\] for all $a\in M_n(\mathcal V)$.

\begin{theorem}\label{th:612} Let $(\mathcal V, 1)$ be a complete matrix order unit space,
and let ${\mathcal D}= (D_n)$ be a matrix metric on  ${\mathcal
C\mathcal S}(\mathcal V)$ such that the $\mathcal{D}$-topology on
$\mathcal{CS(V)}$ agrees with the BW-topology. Then $\mathcal D$
comes from a lower semicontinuous matrix Lip-gauge
$\mathcal{G}=(G_n)$ on $(\mathcal V, 1)$, via the relation
\[D_n(\varphi, \psi)=D_{G_n}(\varphi,  \psi),\,\hbox{for }\varphi,\psi\in CS_n(\mathcal{V}),\,
n\in\mathbb{N},\] if and only if $\mathcal D$ is convex, midpoint
balanced, and midpoint concave.
\end{theorem}

\begin{proof} If $\mathcal D$ comes from a lower semicontinuous matrix
Lip-gauge $\mathcal{G}=(G_n)$ on $(\mathcal V, 1)$ via the
relation $D_n(\varphi, \psi)=D_{G_n}(\varphi, \psi)$ for
$\varphi,\psi\in CS_n(\mathcal{V})$, it is clear that $\mathcal D$
is convex, midpoint balanced, and midpoint concave. And the remain follows from Theorem \ref{th:610}.
\end{proof}

\section{The pre-dual of $(\widetilde{\mathcal V}, \widetilde{\mathcal L})$}\label{sec:7}

Suppose $(\mathcal V, 1)$ is a matrix order unit space with operator space norm $\|\cdot\|=
(\|\cdot\|_n)$. Let $(\overline{\mathcal V}, 1)$ denote the completion of $(\mathcal V, 1)$
for $\|\cdot\|_1$. Given a matrix Lipschitz seminorm $\mathcal L=(L_n)$ on $(\mathcal V, 1)$,
we will denote by $\widehat{\mathcal L^1}$ the norm closure of $\mathcal L^1$ in
$(\overline{\mathcal V}, 1)$. Since in the operator space $\mathcal V$, a sequence
$a^{(k)}=[a^{(k)}_{ij}]$ in $M_n(\mathcal V)$ converges if and only if the entries
$a^{(k)}_{ij}$ converge, we have that for all $m, n\in\mathbb N$, $\widehat{L^1_m}
\oplus\widehat{L^1_n}\subseteq\widehat{L^1_{m+n}}$ and $\alpha\widehat{L^1_m}\beta
\subseteq\widehat{L^1_n}$ for any contractions $\alpha\in M_{n, m}$ and $\beta\in M_{m, n}$.
Thus $\widehat{\mathcal L^1}$ is absolutely matrix convex. So the corresponding Minkowski
functionals
\[L_n^-(v)=\inf\left\{\lambda>0:\, v\in\lambda\widehat{L_n^1}\right\}\]
associated with the convex sets $\widehat{L_n^1}$ determine a matrix gauge $\mathcal L^-=
(L_n^-)$ on $(\overline{\mathcal V}, 1)$(see page 171 in \cite{efwe}). $\widehat{\mathcal L^1}$ is closed
implies that $\mathcal L^-$ is lower semicontinuous.

\begin{proposition}\label{pro:71}
Let $\mathcal L=(L_n)$ be a matrix Lipschitz seminorm on a matrix
order unit space $(\mathcal V, 1)$. Then
\begin{enumerate}
\item $\mathcal{L}$ is lower semicontinuous if and only if
$\mathcal L^-$ is an extension of $\mathcal L$;
\item $\mathcal{D}_{\mathcal{L}^-}=\mathcal{D_L}$ on
$\mathcal{CS}(\overline{\mathcal{V}})=\mathcal{CS}(\mathcal{V})$.
\end{enumerate}
\end{proposition}

\begin{proof} (1) Suppose that $\mathcal{L}$ is lower semicontinuous and $a\in M_n(\mathcal V)$.
If $L_n(a)\le 1$, then $a\in L^1_n\subseteq\widehat{L^1_n}$, and
so clearly $L_n^-(a)\le 1$. If $L_n^-(a)\le 1$, then
$a\in\widehat{L_n^1}$. Thus there is a sequence $\{a_k\}$ in
$L_n^1$ which converges to $a$. From the lower semicontinuity of
$L_n$ it follows that $L_n(a)\le 1$. So $L_n(a)=L_n^-(a)$ for
$a\in M_n(\mathcal V)$.

Conversely, if $\mathcal L^-$ is an extension of $\mathcal L$,
then for any $t>0$, $L_n^t=\{a\in M_n(\mathcal{V}): L_n(a)\le
t\}=\{a\in M_n(\mathcal{V}): L_n^-(a)\le t\}=L_n^{-t}\cap
M_n(\mathcal{V})$, and hence $L_n^t$ is closed in
$M_n(\mathcal{V})$.

(2) If $\varphi,\psi\in CS_n(\mathcal{V})$, we have
\[\begin{array}{rcl}
D_{L_n^-}(\varphi,\psi)
&=&\sup\{\|\ll\varphi,a\gg-\ll\psi,a\gg\|: a\in\widehat{L^1_r}, r\in\mathbb{N}\}\\
&=&\sup\{\|\ll\varphi,a\gg-\ll\psi,a\gg\|: a\in L^1_r,
r\in\mathbb{N}\}\\
&=&D_{L_n}(\varphi,\psi),\end{array}\] since
$CS_n(\overline{\mathcal{V}})=CS_n(\mathcal{V})$.
\end{proof}

\begin{definition}\label{def:72} Let $(\mathcal V, 1)$ be a matrix order unit space, and let $\mathcal L=
(L_n)$ be a matrix Lipschitz seminorm on $(\mathcal V, 1)$. The matrix gauge $\mathcal L^-=
(L_n^-)$ on $(\overline{\mathcal V}, 1)$ is called the {\it closure} of $\mathcal L$. We say
that a matrix Lipschitz seminorm $\mathcal L$ is {\it closed} if $\mathcal L=\mathcal L^-$ on
the subspace where $\mathcal L^-$ is finite.
\end{definition}

\begin{lemma}\label{le:73} For $b\in M_n(\bar{\mathcal V})$, $L_n^-(b)\le 1$ if and only if $b\in
\widehat{L_n^1}$.\end{lemma}

\begin{proof} It is clear that $L_n^-(b)\le 1$ if $b\in\widehat{L_n^1}$. Now suppose that
$L_n^-(b)\le 1$. Then there is a sequence $\{\lambda_k\}$ in
$\mathbb R$ such that $L_n^-(b)\le\lambda_k\le L_n^-(b)+\frac1k$
and $b\in\lambda_k\widehat{L_n^1}$. Let $b=\lambda_kb_k$ for some
$b_k\in\widehat{L_n^1}$ and let $c_k\in L_n^1$ such that
$\|b_k-c_k\|_n\le\frac1k$. Then
\[\|b-\lambda_kc_k\|_n=\lambda_k\|b_k-c_k\|_n\le\frac1k\left(L_n^-(b)+\frac1k\right)\to 0,\ \ \ \
(k\to +\infty).\] If $L_n^{-}(b)=0$, then
$L_n(\lambda_kc_k)=\lambda_kL_n(c_k)\le k^{-1}\le 1$, and hence
$b=\lim_{k\to\infty}\lambda_kc_k\in \widehat{L_n^1}$. If
$L_n^{-}(b)>0$, then
\[\begin{array}{rcl}
\|b-L_n^-(b)c_k\|_n&\le&\|b-\lambda_kc_k\|_n+\|(\lambda_k-L_n^-(b))c_k\|_n\\
&=&\lambda_k\|b_k-c_k\|_n+(1-L_n^-(b)\lambda_k^{-1})\|\lambda_kc_k\|_n\\
&\le&\frac1k(L_n^-(b)+\frac1k)+(kL_n^{-}(b))^{-1}\|\lambda_kc_k\|_n\\
&\to& 0\ \ \ \ (k\to +\infty).\end{array}\] So $b=\lim_{k\to\infty}
L_n^-(b)c_k$. Since $L_n(L_n^-(b)c_k)=L_n^-(b)L_n(c_k)\le 1$, we
have that $b\in \widehat{L_n^1}$.
\end{proof}

\begin{lemma}\label{le:74} $\mathcal{L}$ is closed if and only if $\widehat{L_n^1}=L_n^1$.\end{lemma}

\begin{proof} If $\widehat{L_n^1}=L_n^1$, then for $b\in M_n(\overline{\mathcal V})$ with
$L_n^-(b)<+\infty$, we have that $b\in M_n(\mathcal V)$ and
$L_n^-(b)=\inf\left\{r\in\mathbb R^+:\, b\in
r\widehat{L_n^1}\right\}=\inf\{r\in\mathbb R^+:\, b\in
rL_n^1\}=L_n(b)$. Conversely, assume that $\mathcal{L}$ is closed.
Then for $b\in\widehat{L_n^1}$ we have that $L_n^-(b)\le 1$ by
Lemma \ref{le:73}. Hence $b\in M_n(\mathcal V)$ and
$L_n(b)=L_n^-(b)\le 1$. So $b\in L_n^1$. Therefore,
$\widehat{L_n^1}=L_n^1$.
\end{proof}

Suppose $(\mathcal V, 1)$ is a matrix order unit space with
operator space norm $\|\cdot\|=(\|\cdot\|_n)$ and
$\mathcal{L}=(L_n)$ be a matrix Lipschitz seminorm on it. Denote
$B^t_n=\{a\in M_n(\mathcal V):\, L_n(a)\le 1, \|a\|_n\le t\}$ for
$t\ge 0$ and $n\in\mathbb N$. We say that $\mathbf B^t=(B^t_n)$ is
totally bounded for $\|\cdot\|$ if each $B^t_n$ is totally bounded
for $\|\cdot\|_n$. $\mathcal{D_L}= (D_{L_n})$ is bounded if every
$D_{L_n}$ is bounded. Clearly, $\mathbf B^t$ is totally bounded
for $\|\cdot\|$ if and only if $B^t_1$ is totally bounded for
$\|\cdot\|_1$.

\begin{proposition}\label{pro:75}
Let $(\mathcal V, 1)$ be a matrix order unit space with operator
space norm $\|\cdot\|=(\|\cdot\|_n)$ and $\mathcal{L}=(L_n)$ be a
matrix Lipschitz seminorm on it. Then the following conditions are
equivalent:
\begin{enumerate}
\item The $\mathcal{D_L}$-topology on $\mathcal{CS(V)}$ agrees
with the BW-topology.
\item The image of $L_1^1$ in $\tilde{\mathcal{V}}$ is totally
bounded for $\|\cdot\|_1^\sim$.
\item $D_{L_1}$ is bounded and $B_1^1$ is totally bounded in
$\mathcal{V}$ for $\|\cdot\|_1$.
\end{enumerate}
\end{proposition}

\begin{proof}By Theorem 5.3 in \cite{wu}, (1) and (2) are
equivalent. The equivalence of (3) and (2) follows from Theorem
1.8 and Theorem 1.9 in \cite{ri0} and Proposition 3.1 in
\cite{wu}.
\end{proof}

Let $(\mathcal V, 1)$ be a matrix order unit space, and let
$\mathcal{L}=(L_n)$ be a matrix Lip-norm on it which is closed.
Denote \[K_n=\{\tilde{a}\in M_n(\tilde{\mathcal{V}}):
\tilde{L}_n(\tilde{a})\le 1\},\ \ \ n\in\mathbb{N}.\] Then
$\mathbf K=(K_n)$ is an absolutely matrix convex set in
$\tilde{\mathcal{V}}$. Because $\mathcal{L}$ is closed, the
totally bounded set $B_n^t$ are complete for $\|\cdot\|_n$ by
Lemma \ref{le:74} and Proposition \ref{pro:75}, so are compact. By
Proposition 3.8 and Proposition 5.2 in \cite{wu}, there is a
sequence $\{r_n\}$ of positive constants such that
\[\|\tilde{a}\|_n^\sim\le r_n\tilde{L}_n(\tilde{a}),\ \ \ \hbox{ for }a\in M_n(\mathcal{V})
\hbox{ and }n\in\mathbb{N}.\] Now we see that there is a $t_n>0$
such that \[K_n=\{\tilde{a}\in M_n(\tilde{\mathcal{V}}): L_n(a)\le
1, \|a\|_n\le t_n\}=\{\tilde{a}\in M_n(\tilde{\mathcal{V}}): a\in
B_n^{t_n}\}.\] Since the quotient mapping $\pi_n:
(M_n(\mathcal{V}),
\|\cdot\|_n)\longmapsto(M_n(\tilde{\mathcal{V}}),\|\cdot\|_n^\sim)$
is a contraction, the image $\pi_n(B_n^t)$ of the compact set
$B_n^t$ is compact. So $K_n$ is compact for $\|\cdot\|_n^\sim$.

Set \[G_n=\{\tilde{a}\in M_n(\tilde{\mathcal{V}}):
\tilde{a}^*=\tilde{a},\tilde{L}_n(\tilde{a})\le 1\}, \ \
n\in\mathbb{N}.\] Then $\mathbf{G}=(G_n)$ is a matrix convex set
in $\tilde{\mathcal V}$ and $\mathbf{G}\subseteq\mathbf{K}$. The
involution in $\tilde{\mathcal{V}}$ preserves the matrix norm
$\|\cdot\|^\sim$ implies that $\mathbf{G}$ is also compact for
$\|\cdot\|^\sim$. For $\tilde{a}\in K_n$, we have
\[\tilde{a}=\left[\begin{array}{cc}
1&0\end{array}\right]\left[\begin{array}{cc}
0&\tilde{a}\\
\tilde{a}^*&0\end{array}\right]\left[\begin{array}{c}
0\\
1\end{array}\right].\] So by Lemma 3.2 in \cite{efwe},
$\mathbf{K}$ is the smallest absolutely matrix convex set
containing $\mathbf{G}$, and
\[\begin{array}{rcl}
K_n&=&\Big\{\sum_{j=1}^{m}\alpha_jg_j\beta_j: \alpha_j\in
M_{n,n_j},
g_j\in G_{n_j}, \beta\in M_{n_j,n},\\
&&\sum_{j=1}^{m}\alpha_j\alpha_j^*\le 1_n,
\sum_{j=1}^m\beta^*_j\beta_j\le 1_n,
m\in\mathbb{N}\Big\}.\end{array}\]

Denote
\[A_0(\mathbf{G})=\{\mathbf{F}=(F_n)\in A(\mathbf{G}): F_n(\tilde{0}_n)=0_n,\ \ n\in\mathbb{N}\}.\]
For $\mathbf{F}^{(k)}=(F_n^{(k)})\in M_k(A_0(\mathbf{G}))\cong
A_0(\mathbf{G}, M_k)$, we define
$\tilde{\mathbf{F}}^{(k)}=(\tilde{F}_n^{(k)})$ on $\mathbf{K}$ by
\[\tilde{F}_n^{(k)}\left(\sum_{j=1}^m\alpha_jg_j\beta_j\right)=\sum_{j=1}^m(\alpha_j\otimes 1_k)F_{n_j}^{(k)}(g_j)
(\beta_j\otimes 1_k),\] where $\alpha_j\in M_{n,n_j}$, $g_j\in
G_{n_j}$, $\beta_j\in M_{n_j,n}$,
$\sum_{j=1}^m\alpha_j\alpha_j^*\le 1_n$ and
$\sum_{j=1}^m\beta_j^*\beta_j\le 1_n$. If
$\sum_{j=1}^m\alpha_jg_j\beta_j=\tilde{0}_n$, then
\[\begin{array}{rcl}
&&\sum_{j=1}^m\frac{1}{2}\left((\alpha_j^*+\beta_j)^*g_j(\alpha_j^*+\beta_j)-\alpha_jg_j\alpha_j^*
-\beta_j^*g_j\beta_j\right)\\
&&-i\sum_{j=1}^m\frac{1}{2}\left((\alpha_j^*+i\beta_j)^*g_j(\alpha_j^*+i\beta_j)
-\alpha_jg_j\alpha_j^*-\beta_j^*g_j\beta_j\right)=\tilde{0}_n,\end{array}\]
that is,
\[\sum_{j=1}^m\frac{1}{2}\left((\alpha_j^*+\beta_j)^*g_j(\alpha_j^*+\beta_j)-\alpha_jg_j\alpha_j^*
-\beta_j^*g_j\beta_j\right)=\tilde{0}_n,\]
\[\sum_{j=1}^m\frac{1}{2}\left((\alpha_j^*+i\beta_j)^*g_j(\alpha_j^*+i\beta_j)
-\alpha_jg_j\alpha_j^*-\beta_j^*g_j\beta_j\right)=\tilde{0}_n.\]
So
\[\sum_{j=1}^m\frac{1}{2(a+1)}(\alpha_j^*+\beta_j)^*g_j(\alpha_j^*+\beta_j)=\sum_{j=1}^m\left(
\frac{1}{2(a+1)}\alpha_jg_j\alpha_j^*+\frac{1}{2(a+1)}\beta_j^*g_j\beta_j\right),\]
\[\sum_{j=1}^m\frac{1}{2(b+1)}(\alpha_j^*+i\beta_j)^*g_j(\alpha_j^*+i\beta_j)=\sum_{j=1}^m\left(\frac{1}{2(b+1)}
\alpha_jg_j\alpha_j^*+\frac{1}{2(b+1)}\beta_j^*g_j\beta_j\right),\]
where
$a=\|\sum_{j=1}^m\frac{1}{2}(\alpha_j^*+\beta_j)^*(\alpha_j^*+\beta_j)\|$
and
$b=\|\sum_{j=1}^m\frac{1}{2}(\alpha_j^*+i\beta_j)^*(\alpha_j^*+i\beta_j)\|$.
Since $\tilde{0}_n\in G_n$ and $\mathbf{F}^{(k)}\in
A_0(\mathbf{G}, M_k)$, we obtain
\[\begin{array}{rcl}
&&\sum_{j=1}^m\frac{1}{2(a+1)}\left((\alpha_j^*+\beta_j)^*\otimes
1_k\right)F_{n_j}^{(k)}(g_j)\left((\alpha_j^*+\beta_j)\otimes 1_k\right)\\
&=&\sum_{j=1}^m\Big(\frac{1}{2(a+1)}(\alpha_j\otimes 1_k)
F_{n_j}^{(k)}(g_j)(\alpha_j^*\otimes 1_k)\\
&&+\frac{1}{2(a+1)}(\beta_j^*\otimes
1_k)F_{n_j}^{(k)}(g_j)(\beta_j\otimes 1_k)\Big),\end{array}\]
and
\[\begin{array}{rcl}
&&\sum_{j=1}^m\frac{1}{2(b+1)}\left((\alpha_j^*+i\beta_j)^*\otimes
1_k\right)F_{n_j}^{(k)}(g_j)\left((\alpha_j^*+i\beta_j)\otimes 1_k\right)\\
&=&\sum_{j=1}^m\Big(\frac{1}{2(b+1)}(\alpha_j\otimes 1_k)
F_{n_j}^{(k)}(g_j)(\alpha_j^*\otimes 1_k)\\
&+&\frac{1}{2(b+1)}(\beta_j^*\otimes
1_k)F_{n_j}^{(k)}(g_j)(\beta_j\otimes 1_k)\Big).\end{array}\] From
these, we see that $\sum_{j=1}^m(\alpha_j\otimes
1_k)F_{n_j}^{(k)}(g_j)(\beta_j\otimes 1_k)=\tilde{0}_{nk}$.
Therefore, $\tilde{\mathbf{F}}^{(k)}$ is well-defined.

\begin{proposition}\label{pro:76} Denote $A_0(\mathbf{K})=\{\tilde{\mathbf{F}}:\mathbf{F}\in
A_0(\mathbf{G})\}$. Then for $\tilde{\mathbf{F}}^{(k)}\in
M_k(A_0(\mathbf{K}))$, $k_i\in K_{n_i}$ and $\alpha_i\in
M_{n,n_i}$, $\beta_i\in M_{n_i,n}$ satisfying
$\sum_{i=1}^m\alpha_i\alpha_i^*\le 1_n$ and
$\sum_{i=1}^m\beta_i^*\beta_i\le 1_n$, we have
\[\tilde{F}_n^{(k)}\left(\sum_{i=1}^m\alpha_ik_i\beta_i\right)=\sum_{i=1}^m(\alpha_i\otimes 1_k)\tilde{F}_{n_i}^{(k)}(k_i)(\beta_i\otimes 1_k).\]
In particular, $A_0(\mathbf{K})\subseteq A(\mathbf{K})$.
\end{proposition}

\begin{proof} Direct verification according to
definition.\end{proof}

\begin{proposition}\label{pro:77}
For $\tilde{\mathbf{F}}^{(k)}\in M_k(A_0(\mathbf{K}))$, we define
\[\check{F}_n^{(k)}\left(\sum_{i=1}^m\lambda_ik_i\right)=\sum_{i=1}^m\lambda_i\tilde{F}_n^{(k)}(k_i)\]
for $\lambda_i\in\mathbb{R}^{+}$ and $k_i\in K_n$,
$i=1,2,\cdots,m$. Then:
\begin{enumerate}
\item $\check{F}_n^{(k)}(\sum_{i=1}^m\alpha_ik_i\beta_i)=\sum_{i=1}^m(\alpha_i\otimes
1_k)\check{F}_{n_i}^{(k)}(k_i)(\beta_i\otimes 1_k)$ for $k_i\in
K_{n_i}$ and $\alpha_i\in M_{n,n_i}$, $\beta_i\in M_{n_i,n}$;
\item If we denote $L_0(\mathbf{K})=\{\check{\mathbf{F}}: \tilde{\mathbf{F}}\in
A_0(\mathbf{K})\}$and still view $L_0(\mathbf{K})$ as equipped
with the operator space norm
$\|\cdot\|^{\sharp}=(\|\cdot\|_n^{\sharp})$ determined by the
matrix order structure on $A(\mathbf{K})$, then
\[\|\check{\mathbf{F}}^{(r)}\|_r^{\sharp}=\sup\{\|\tilde{F}_n^{(r)}(\tilde{a})\|:
\tilde{a}\in K_n
, n\in\mathbb{N}\},\] for $\mathbf{F}^{(r)}=(F_n^{(r)})\in
A_0(\mathbf{G}, M_r )$ and $r\in\mathbb{N}$.
\end{enumerate}
\end{proposition}

\begin{proof}
Suppose that $\sum_{i=1}^m\lambda_ik_i=\sum_{j=1}^l\mu_jh_j$ with
$k_i, h_j\in K_n$ and $\lambda_i, \mu_j\in\mathbb{R}^{+}$,
$i=1,2,\cdots,m$, $j=1,2,\cdots,l$. Let
$\sum_{i=1}^m\lambda_i+\sum_{j=1}^l\mu_j=a>0$. Then
$\sum_{i=1}^m\frac{\lambda_i}{a}k_i=\sum_{j=1}^l\frac{\mu_j}{a}h_j$.
Since $\tilde{0}_n\in K_n$,
$\sum_{i=1}^m\frac{\lambda_i}{a}\tilde{F}_n^{(k)}(k_i)=\sum_{j=1}^l\frac{\mu_j}{a}\tilde{F}_n^{(k)}(h_j)$
by Proposition \ref{pro:76}, that is,
$\sum_{i=1}^m\lambda_i\tilde{F}_n^{(k)}(k_i)=\sum_{j=1}^l\mu_j\tilde{F}_n^{(k)}(h_j)$.
Thus $\check{\mathbf{F}}^{(k)}=(\check{F}_n^{(k)})$ is
well-defined.

For $k_i\in K_{n_i}$ and $\alpha_i\in M_{n,n_i}$, $\beta_i\in
M_{n_i,n}$, $i=1,2,\cdots, m$, we denote
\[\max\left\{\left\|\sum_{i=1}^m\alpha_i\alpha_i^*\right\|, \ \ \ \ \left\|\sum_{i=1}^m\beta_i^*\beta_i\right\|\right\}=b.\]
Then
\[\sum_{i=1}^m\frac{\alpha_i}{\sqrt{b+1}}\left(\frac{\alpha_i}{\sqrt{b+1}}\right)^*\le 1_n, \ \ \
\sum_{i=1}^m\left(\frac{\beta_i}{\sqrt{b+1}}\right)^*\frac{\beta_i}{\sqrt{b+1}}\le
1_n.\] By Proposition \ref{pro:76},
\[\begin{array}{rcl}
\check{F}_n^{(k)}\left(\sum_{i=1}^m\alpha_ik_i\beta_i\right)
&=&\check{F}_n^{(k)}\left((b+1)\sum_{i=1}^m\frac{\alpha_i}{\sqrt{b+1}}k_i\frac{\beta_i}{\sqrt{b+1}}\right)\\
&=&(b+1)\tilde{F}_n^{(k)}\left(\sum_{i=1}^m\frac{\alpha_i}{\sqrt{b+1}}k_i\frac{\beta_i}{\sqrt{b+1}}\right)\\
&=&\sum_{i=1}^m(\alpha_i\otimes 1_k)\tilde{F}_n^{(k)}(k_i)(\beta_i\otimes 1_k)\\
&=&\sum_{i=1}^m(\alpha_i\otimes
1_k)\check{F}_n^{(k)}(k_i)(\beta_i\otimes 1_k).
\end{array}\]
This proves (1).

For $\tilde{\mathbf F}^{(r)}\in A_0(\mathbf K, M_r)$, we have
\[\begin{array}{rcl}
\|\check{\mathbf F}^{(r)}\|_r^\sharp&=& \|\tilde{\mathbf
F}^{(r)}\|_r^\sharp=\inf\left\{t>0:\, \left[\begin{array}{cc}
t\mathbf{I}^{(r)}&\tilde{\mathbf F}^{(r)}\\
\tilde{\mathbf F}^{^{(r)}*}&t{\mathbf I}^{(r)}\end{array}\right]\ge 0\right\}\\
&=&\inf\left\{t>0:\, \left[\begin{array}{cc}
t(1_n\otimes 1_r)&\tilde{F}^{(r)}_n(\tilde{a})\\
\tilde{F}^{(r)}_n(\tilde{a})^*&t(1_n\otimes 1_r)\end{array}\right]\ge 0, \tilde{a}\in K_n, n\in\mathbb N\right\}\\
&=&\sup\{\|\tilde{F}^{(r)}_n(\tilde{a})\|:\, \tilde{a}\in K_n,
n\in\mathbb N\},
\end{array}\]
and (2) follows.
\end{proof}

\begin{theorem}\label{th:78} Let $(\mathcal V, 1)$ be a matrix order unit space, and let $\mathcal L=
(L_n)$ be a matrix Lip-norm on $(\mathcal V, 1)$ which is closed.
Denote
\[K_n=\{\tilde{a}\in M_n(\tilde{\mathcal V}):\,\tilde{L}_n(\tilde{a})\le 1\},\]
and set $\mathbf{K}=(K_n)$. Then $(\tilde{\mathcal V}, \tilde{\mathcal{L}})$
is naturally completely isometrically isomorphic to the dual
operator space of $A_0(\mathbf{K})$.
\end{theorem}

\begin{proof} Let $\|\cdot\|^{\sharp}=(\|\cdot\|_n^{\sharp})$ be the operator space norm on
$A(\mathbf{K})$.Clearly $A_0(\mathbf{K})$ is a self-adjoint
subspace of $A(\mathbf{K})$. Let $A_0(\mathbf{K})^{\uparrow}$ be
the operator space dual of $A_0(\mathbf{K})$ with the dual
operator space norm
$\|\cdot\|^{\uparrow}=(\|\cdot\|_n^{\uparrow})$. For
$\tilde{a}\in\tilde{\mathcal{V}}$, we define
\[\tau (\tilde{a})(\tilde{\mathbf{F}}^{(r)})=\check{F}_1^{(r)}(\tilde{a}),\]
where $\tilde{\mathbf{F}}^{(r)}\in M_r(A_0(\mathbf{K}))$. By
Proposition \ref{pro:77}, $\tau$ is well-defined and
$\tau(\tilde{a})$ is linear. If we choose $t>L_1(\tilde{a})$, then
\[\|\tau (\tilde{a})(\tilde{\mathbf F}^{(r)})\|=\left\|t\tilde{F}_1^{(r)}\left(\frac{\tilde{a}}{t}\right)\right\|
\le t\|\tilde{\mathbf{F}}^{(r)}\|_r^{\sharp},\] for
$\tilde{\mathbf F}^{(r)}\in M_r(A_0(\mathbf K))$. So
$\tau(\tilde{a})\in A_0(\mathbf{K})^{\uparrow}$ and
$\|\tau(\tilde{a})\|_1^{\uparrow}\le L_1(\tilde{a})$. From
Proposition \ref{pro:77}, we see that $\tau$ is linear.

To see injectivity, let $(\tilde{\mathcal{V}})^*$ be the dual
Banach space of $\tilde{\mathcal{V}}$, and then observe that $f\in
(\tilde{\mathcal{V}})^*$ defines an element in $A(\mathbf{K})$
determined by the linear mapping
$\tilde{a}\in\tilde{\mathcal{V}}\longmapsto f(\tilde{a})$. If
$\tau(\tilde{a})=\tau(\tilde{b})$ for $\tilde{a},
\tilde{b}\in\tilde{\mathcal{V}}$, then in particular
$f(\tilde{a})=f(\tilde{b})$ for all $f\in
(\tilde{\mathcal{V}})^{*}$, which again implies that
$\tilde{a}=\tilde{b}$, since $\tilde{\mathcal{V}}^*$ separates
points in $\tilde{\mathcal{V}}$.

For $\tilde{a}=[\tilde{a}_{ij}]\in M_n(\tilde{\mathcal{V}})$ and
$\tilde{\mathbf{F}}^{(r)}\in A_0(\mathbf{K}, M_r)$, we have
\[\begin{array}{rcl}
\tau_n(\tilde{a})(\tilde{\mathbf{F}}^{(r)})&=&\ll\tilde{\mathbf{F}}^{(r)},
[\tau(\tilde{a}_{ij})]\gg=[\ll\tilde{\mathbf{F}}^{(r)},
\tau(\tilde{a}_{ij})\gg]\\
&=&[\check{F}_1^{(r)}(\tilde{a}_{ij})]
=\check{F}_n^{(r)}(\tilde{a})=
(\tilde{L}_n(\tilde{a})+\epsilon)\tilde{F}_n^{(r)}
\left(\frac{\tilde{a}}{\tilde{L}_n(\tilde{a})+\epsilon}\right)\end{array}\]
for $\epsilon>0$. So
$\|\tau_n(\tilde{a})(\tilde{\mathbf{F}}^{(r)})\|\le
\tilde{L}_n(\tilde{a})\|\tilde{\mathbf{F}}^{(r)}\|_r^\sharp$, that
is, $\|\tau_n(\tilde{a})\|_n^\uparrow\le\tilde{L}_n(\tilde{a})$.
Especially, $\tau_n(K_n)\subseteq
M_n(A_0(\mathbf{K})^\uparrow)_{\|\cdot\|_n^\uparrow\le 1}$, where
$M_n(A_0(\mathbf{K})^\uparrow)_{\|\cdot\|_n^\uparrow\le 1}=\{f\in
M_n(A_0(\mathbf{K})^\uparrow): \|f\|_n^\uparrow\le 1\}$. Since for
$k_1, k_2\in K_n$,
\[\|(\tau_n(k_1)-\tau_n(k_2))(\check{\mathbf{F}}^{(r)})\|=\|\check{F}_n^{(r)}(k_1)-\check{F}_n^{(r)}(k_2)\|,\]
$\tau_n$ is continuous from $K_n$ to
$M_n(A_0(\mathbf{K})^\uparrow)_{\|\cdot\|_n^\uparrow\le 1}$ with
the point-norm topology. $K_n$ is compact implies that
$\tau_n(K_n)$ must be compact for the point-norm topology.

Assume that $f_0^*=f_0\in
M_n(A_0(\mathbf{K})^\uparrow)_{\|\cdot\|_n^\uparrow\le
1}\setminus\tau_n(K_n)$. By the matricial separation theorem,
there is a continuous linear mapping $\Phi:
A_0(\mathbf{K})^\uparrow\longmapsto M_n$ such that
\[\mathrm{Re}\Phi_r(\tau_r(k))\le 1_n\otimes1_r\]
for all $k\in K_r$, $r\in\mathbb{N}$, and
\[\mathrm{Re}\Phi_n(f_0)\not\leq 1_n\otimes 1_n.\]
Identifying $\Phi$ with $\mathbf{F}^{(n)}\in
M_n(A_0(\mathbf{K}))\cong A_0(\mathbf{K}, M_n)$, this means that
\[\mathrm{Re}\mathbf{F}_r^{(n)}(k)\le 1_n\otimes 1_r\]
for all $k\in K_r$, $r\in\mathbb{N}$, and
\[\ll f_0, \mathrm{Re}\mathbf{F}^{(n)}\gg=\mathrm{Re}\ll f_0, \mathbf{F}^{(n)}\gg\not\le 1_n\otimes 1_n,\]
since $f_0$ is self-adjiont. It is clear that
$\mathrm{Re}\mathbf{F}^{(n)}\in M_n(A_0(\mathbf{K}))$, and the
first inequality and Proposition \ref{pro:77} say that
$\|\mathrm{Re}\mathbf{F}^{(n)}\|^\sharp_n \le 1$. The second
inequality implies that $\|f_0\|_n^\uparrow>1$, contradicting that
$f_0\in M_n(A_0(\mathbf{K})^\uparrow)_{\|\cdot\|_n^\uparrow\le
1}$. Now suppose that $f_0\in
M_n(A_0(\mathbf{K})^\uparrow)_{\|\cdot\|_n^\uparrow\le 1}$. Then
$\left[\begin{array}{cc}0&f_0\\ f_0^*&0\end{array}\right]$ is
self-adjoint and belongs to
$M_{2n}(A_0(\mathbf{K})^\uparrow)_{\|\cdot\|_{2n}^\uparrow\le 1}$.
So there is an element $\left[\begin{array}{cc}k_{11}&k_{12}\\
k_{21}&k_{22}\end{array}\right]\in K_{2n}$ such that
\[\tau_{2n}\left(\left[\begin{array}{cc}k_{11}&k_{12}\\
k_{21}&k_{22}\end{array}\right]\right)=\left[\begin{array}{cc}0&f_0\\
f_0^*&0\end{array}\right].\] Thus
$f_0=\tau_n(k_{12})\in\tau_n(K_n)$. Therefore,
$\tau_n(K_n)=M_n(A_0(\mathbf{K})^\uparrow)_{\|\cdot\|_n^\uparrow\le
1}$. Consequently $\tau$ is a completely isometric isomorphism of
$(\tilde{\mathcal{V}}, \tilde{\mathcal{L}})$ onto
$(A_0(\mathbf{K})^\uparrow,\|\cdot\|^\uparrow)$.
\end{proof}

\section*{Acknowledgements}

This research was partially supported by Shanghai Priority Academic Discipline, China Scholarship Council, and National Natural Science 
Foundation of China. I would like to thank Marc Rieffel for many helpful discussions and suggestions. I want to thank Hanfeng Li for 
valuable suggestions and comments.

\bibliographystyle{amsplain}

\end{document}